\documentclass[a4paper,10pt]{article}
\usepackage{kubik}
\newcommand{\nofootnote}[1]{%
\begingroup\def\thefootnote{}\footnotetext{#1}\endgroup}

\renewcommand*{\HyperDestNameFilter}[1]{\jobname-#1}

\begin{document}
\title{On the total curvature of \\
minimizing geodesics on convex surfaces}
\author{Nina Lebedeva and Anton Petrunin}
\nofootnote{1991 Mathematics Subject Classification. 53A05, 53C45, 53C22.}
\nofootnote{N.~Lebedeva was partially supported by RFBR grant 14-01-00062.}
\nofootnote{A.~Petrunin was partially supported by NSF grant DMS 1309340.}

\newcommand{\Addresses}{{\bigskip\footnotesize

  Nina Lebedeva, \par\nopagebreak\textsc{Math. Dept.
St. Petersburg State University,
Universitetsky pr., 28, 
Stary Peterhof, 
198504, Russia.}
  \par\nopagebreak
  \textsc{Steklov Institute,
27 Fontanka, St. Petersburg, 
191023, Russia.}
  \par\nopagebreak
  \textit{Email}: \texttt{lebed@pdmi.ras.ru}

\medskip

  Anton Petrunin, \par\nopagebreak\textsc{Math. Dept. PSU, University Park, PA 16802, USA}
  \par\nopagebreak
  \textit{Email}: \texttt{petrunin@math.psu.edu}
}}
\date{}
\maketitle

\begin{abstract}
We give a universal upper bound 
for the total curvature 
of minimizing geodesic 
on a convex surface 
in the Euclidean space.
\end{abstract}

\section{Introduction}

In this note we give an affirmative answer to the question asked by Dmitry Burago; 
later we found the same question in \cite{AH-PSV}, \cite{pach} and \cite{BKZ}.
Namely, we prove the following.

\begin{thm}{Main theorem}\label{thm:main}
The total curvature of a minimizing geodesic
on a convex surface in $\RR^3$ can not exceed $1000^{1000}$.
\end{thm}

The value $2\cdot\pi$ is the optimal bound for the analogous problem in the plane.
The total curvature
of minimizing geodesic on a convex surface in $\RR^3$
can exceed $2\cdot\pi$
and the optimal bound 
is expected to be slightly bigger than $2\cdot\pi$.
The former example was constructed by B{\'a}r{\'a}ny,
Kuperberg, 
and Zamfirescu in \cite{BKZ}. 

Let us list other related results.

\begin{itemize}
\item In \cite{liberman}, Liberman gives a bound on the total curvature of short geodesic on the surface of convex body $K$
in terms of the ratio diameter and inradius of $K$.
In the proof he use an analog of Lemma \ref{lem:liberman} discussed below.

\item In \cite{usov}, 
Usov gives the optimal bound for the total curvature of geodesics on the graphs of $\ell$-Lipschitz convex function. 
Namely, he proves that if $f\:\RR^2\to\RR$ is $\ell$-Lipschitz and convex, then any 
geodesic in its graph 
\[\Gamma_f=\set{(x,y,z)\in \RR^3}{z=f(x,y)}\] 
has total curvature at most $2\cdot \ell$.
An amusing generalization of Usov's result is given by Berg in \cite{berg}.
\item In \cite{pogorelov}, Pogorelov conjectured 
that the spherical image of geodesic on convex surface has to be (locally) rectifiable.
It is easy to see that the length of spherical image of geodesic can not be smaller than its total curvature, 
so this conjecture (if it would be true) 
would be stronger than Liberman's theorem.
Counterexamples to the different forms of this conjecture were found 
by Zalgaller in \cite{zalgaller},
Milka in \cite{milka}
and Usov in \cite{usov-conj-pog};
these results were partly rediscovered later 
by Pach in \cite{pach}.
\item In \cite{BKZ},
B{\'a}r{\'a}ny,
Kuperberg 
and Zamfirescu 
have constructed a corkscrew minimizing geodesic on a closed convex surface;
that is a minimizing geodesic which twists around a given line arbitrary many times.
They also rediscovered the results of Liberman and Usov mentioned above.
\end{itemize}

\parbf{Idea of the proof.}
First we show that it is sufficient to estimate total curvature for the minimizing geodesics with almost constant velocity vector, say $\dot\gamma(t)\approx\bm{i}$.

To estimate total curvature in this case it is sufficient 
to estimate the integral
\[\int\langle\ddot\gamma(t),\bm{j}\rangle\cdot dt\] 
for a vector $\bm{j}\perp\bm{i}$.
To understand the idea of this estimate, 
imagine that 
the surface is lighten in the direction of $\bm{j}$,
so it is divided into the dark and bright sides 
by a curve $\omega$.
On the diagram you see different combinatorics in which
$\gamma$ meets $\omega$.

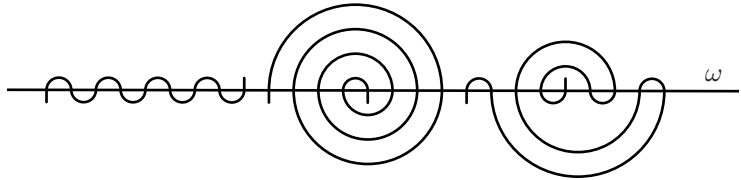
\begin{center}
 \begin{lpic}[t(2 mm),b(2 mm),r(0 mm),l(0 mm)]{pics/spiral-and-snake(1)}
\lbl[b]{93,14;$\omega$}
\end{lpic}
\end{center}

In the first case the total curvature is estimated by integral of Gauss curvature of the regions squeezed between $\gamma$ and $\omega$.
The latter follows from the Tongue Lemma \ref{lem:tongue},
which is the hart of our proof.

The second case might look impossible, 
but the corkscrew geodesic constructed in \cite{BKZ} can meets $\omega$ in this order.
Here we calculate the total curvature for each twists
and show that the obtained sequence of numbers grows geometrically from some place in the middle to the ends. 
At the ends the integral of the full twists can not be bigger than $2\cdot\pi$ ---
this is sufficient to estimate the total curvature of the whole geodesic.

The last case is a mixture of first two and it is done by mixing both techniques.

\section{Preliminaries}

\parbf{Total curvature.}
Recall that the \emph{total curvature} of a curve $\gamma\:[0,\ell]\to \RR^3$ 
(briefly $\tc\gamma$)
is defined as supremum of sum of exterior angles 
for the polygonal lines inscribed in $\gamma$.

Note that for a polygonal line $\sigma$, 
its total curvature coincides with the sum of its exterior angles.

If $\gamma$ is a smooth curve with unit-speed parametrization, 
then 
\[\tc\gamma=\int\limits_0^\ell \kappa(t)\cdot dt,\]
where $\kappa(t)=|\ddot\gamma(t)|$ is the curvature of $\gamma$ at $t$.

\begin{thm}{Proposition}\label{prop:semicontinuity}
Assume $\gamma_n\:\II\to\RR^3$ is a sequence of curves  converging pointwise to a curve $\gamma_\infty\:\II\to\RR^3$.
Then 
\[\liminf_{n\to\infty}\tc\gamma_n\ge \tc\gamma_\infty.\]
\end{thm}

\parit{Proof.}
Fix a polygonal line $\sigma_\infty$ inscribed in $\gamma_\infty$.
Let $\gamma_\infty(t_0),\zz\dots,\gamma_\infty(t_k)$
be its vertices 
as they appear on $\gamma_\infty$.
Consider the polygonal line $\sigma_n$ inscribed in $\gamma_n$ with the vertices 
$\gamma_n(t_0),\zz\dots,\gamma_n(t_k)$.
Note that 
\[\tc\sigma_n\to\tc\sigma_\infty\ \ \text{as}\ \ n\to\infty.\]
By the definition of total curvature, 
\[\tc\sigma_n\le\tc\gamma_n.\]
The statement follows since the broken line $\sigma_\infty$ can be chosen in such a way that 
$\tc\sigma_\infty$ is arbitrary close to $\tc\gamma_\infty$.
\qeds

\parbf{Convergence of sets.}
Given a closed set $\Sigma\subset \RR^3$,
denote by $\dist_\Sigma$ the distance function from $\Sigma$;
that is 
\[\dist_\Sigma x=\inf\set{|x-y|}{y\in\Sigma}.\]

We say that a sequence of closed sets $\Sigma_n\subset \RR^3$
\emph{converges} to the closed set $\Sigma_\infty\subset \RR^3$ 
if for any $x\in\RR^3$,
we have
$\dist_{\Sigma_n}x\to\dist_{\Sigma_\infty}x$ as $n\to\infty$.

\parbf{Convex surfaces.}
By \emph{convex surface} $\Sigma$ in the Euclidean 3-space $\RR^3$ we understand the boundary of a closed convex set $K$ with nonempty interior.
If $K$ is compact we say that the $\Sigma$ is \emph{closed}.

Assume $\Sigma$ is smooth.
If at every point both of $\Sigma$
the principle curvatures are positive, 
we say that $\Sigma$ is \emph{strongly convex}. 

\begin{thm}{Proposition}\label{prop:convegence}
Assume $\Sigma_n$ be a sequence of convex surfaces which converge to a convex surface $\Sigma_\infty$.
Then for any minimizing geodesic $\gamma_\infty$ in $\Sigma_\infty$ there is a sequence of minimizing geodesics $\gamma_n$ in $\Sigma_n$ which pointwise converge to $\gamma_\infty$ as $n\to \infty$.
\end{thm}

\parit{Proof.}
Assume $\gamma_\infty\:[0,\ell]\to\Sigma_\infty$ is parametrized by its arc-length.

Fix a sub-interval $[a,b]\subset (0,\ell)$.
Set $p_\infty=\gamma_\infty(a)$ and $q_\infty=\gamma_\infty(b)$.
Let $p_n,q_n\zz\in\Sigma_n$ be a two sequences of points which converge to $p_\infty$ and $q_\infty$ 
correspondingly.

Denote by $\gamma_n$ a minimizing geodesic from $p_n$ to $q_n$ in $\Sigma_n$.
Note that $\gamma_n$ converges to $\gamma_\infty|_{[a,b]}$
as $n\to\infty$.

Taking the sub-interval  $[a,b]$ closer and closer to $[0,\ell]$
and applying diagonal procedure, we get the result.
\qeds

\section{Liberman's lemma}

In this section we give a slight generalization 
of the construction given by Liberman in \cite{liberman};
see also \cite{milka-liberman} and \cite{petrunin}.

\parbf{Development.}
Let $\gamma\:[0,\ell]\to\RR^3$ be a curve parametrized by arc-length
and a point $p$ does not lie on $\gamma$.

Assume  $\tilde\gamma_p\:[0,\ell]\to\RR^2$ is a plane curve parametrized by arc-length
and $\tilde p$ is a point in the plane such that 
\[|\tilde p-\tilde\gamma(t)|=|p-\gamma(t)|\]
for any $t\in[0,\ell]$.
Moreover, 
the direction from $\tilde p$ to $\tilde \gamma(t)$ changes monotonically (clockwise or counterclockwise). 
Then $\tilde\gamma_p$ is called \emph{development} of $\gamma$ with respect to $p$.

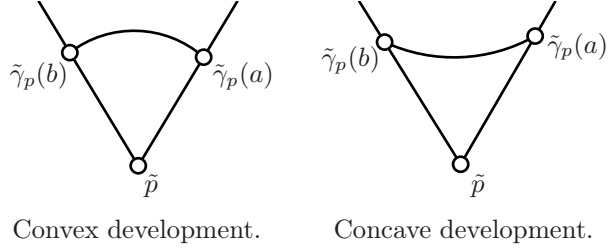
\begin{center}
\begin{lpic}[t(0 mm),b(9 mm),r(0 mm),l(0 mm)]{pics/convex-concave(1)}
\lbl[tl]{23,16;$\tilde\gamma_p(a)$}
\lbl[tr]{4,16;$\tilde\gamma_p(b)$}
\lbl[tl]{14,1;$\tilde p$}
\lbl[ll]{67,18;$\tilde\gamma_p(a)$}
\lbl[rt]{45,18;$\tilde\gamma_p(b)$}
\lbl[tl]{57,1;$\tilde p$}
\lbl[t]{13,-5;Convex development.}
\lbl[t]{56,-5;Concave development.}
\end{lpic}
\end{center}

We say that the development $\tilde\gamma_p$ is \emph{convex} (\emph{concave}) in the interval $[a,b]$
if the arc $\tilde\gamma_p|_{[a,b]}$
cuts from the solid angle $\angle \tilde p^{\tilde\gamma_p(a)}_{\tilde\gamma_p(b)}$
convex bounded (correspondingly unbounded) domain.

We say that $\tilde\gamma_p$ is \emph{locally} convex (concave) in the interval $[a,b]$
if any point $x\in [a,b]$ admits a closed neighborhood $[a',b']$ in $[a,b]$
such that $\tilde\gamma_p$ is convex (correspondingly concave) in the interval $[a',b']$.

If we pass to the limit of this construction as $p$ moves to infinity along a half-line in the direction of a unit vector $-\bm{u}$, then the limit curve is called development of $\gamma$ in the direction $\bm{u}$ and denoted as $\tilde\gamma_{\bm{u}}$.

We can define the development $\tilde\gamma_{\bm{u}}$ directly:
(1) the development $\tilde\gamma_{\bm{u}}\:[0,\ell]\to \R^2$
is parametrized by arc-length,
(2) for a fixed unit vector $\tilde{\bm{u}}\in\RR^2$,
we have
\[\langle \tilde{\bm{u}},\tilde\gamma_{\bm{u}}(t)\rangle
=
\langle  \bm{u},\gamma(t)\rangle\]
for any $t\in [0,\ell]$
and 
(3) the projection of $\tilde\gamma_{\bm{u}}(t)$ to the line normal to $\tilde{\bm{u}}$
is monotonic in $t$.

We can assume that $\tilde{\bm{u}}$ is the vertical vector in the coordinate plane.
In this case we say that $\tilde\gamma_{\bm{u}}$ is concave (convex) in the interval $[a,b]$ 
if the lune bounded by arc $\tilde\gamma_{\bm{u}}|_{[a,b]}$
and the segment $[\tilde\gamma_{\bm{u}}(a)\tilde\gamma_{\bm{u}}(b)]$
is convex and lies above (correspondingly below) the line segment $[\tilde\gamma_{\bm{u}}(a)\tilde\gamma_{\bm{u}}(b)]$.
 
\parbf{Dark and bright sides.}
Let $\Sigma\subset\RR^3$ be a convex surface,
$p\in\Sigma$ and $z\ne p$.

We say that $p$ lies on the \emph{dark} (\emph{bright}) side of $\Sigma$ with from $z$ 
if non of the points $p+t\cdot(p-z)$ lie inside of $\Sigma$ for $t>0$ (correspondingly for $t<0$).
The intersection of dark and bright side is called horizon of $z$;
it is denoted by $\omega_z$.

Note that if $z$ lies inside $\Sigma$, 
then all the points on $\Sigma$ lie on the dark side from $z$ and its horizon $\omega_z$ is empty.

If $\Sigma$ is smooth we can define the outer normal vector $\nu_p$ to $\Sigma$ at $p$.
In this case $p$ lies on dark (bright) side of $\Sigma$ from $z$
if and only if $\langle p-z,\nu_p\rangle\ge 0$
(correspondingly $\langle p-z,\nu_p\rangle\le 0$).
If in addition $\Sigma$ is closed and strongly convex, then the horizon is empty for $z$ inside $\Sigma$ 
and it is formed by a closed smooth curve for $z$ outside $\Sigma$.

We could also define bright/dark sides and horizon in the limit case,
as $p$ escapes to infinity along a half-line in direction $-\bm{u}$.

The latter can be defined directly.
We say that a point $p\in\Sigma$ lies on dark (bright) side for the unit vector $\bm{u}$ if non of the points $p+\bm{u}\cdot t$ lie inside of $\Sigma$ for $t>0$, (correspondingly $t<0$).
As before the intersection of bright and dark side is called horizon of $\bm{u}$ and it is denoted as $\omega_{\bm{u}}$.

In the smooth case the latter means that $\langle \nu_p,\bm{u}\rangle\ge 0$ (correspondingly $\langle \nu_p,\bm{u}\rangle\le 0$).
If in addition $\Sigma$ is closed strongly convex, then $\omega_{\bm{u}}$ is a closed smooth curve.

\begin{thm}{Liberman's Lemma}\label{lem:liberman}
Assume $\gamma$ be a geodesic on the convex surface $\Sigma\subset \RR^3$.
Then for any point $z\notin\Sigma$ the development $\tilde\gamma_z$ is locally convex (locally concave) if $\gamma$ lies on dark (correspondingly bright) side of $\Sigma$ from $z$.

Similarly for any unit vector $\bm{u}$,
the development $\tilde\gamma_{\bm{u}}$ is locally convex (locally concave) if $\gamma$ lies on dark (correspondingly bright) side of $\Sigma$ for $\bm{u}$.

\end{thm}

Note that for any space curve $\gamma$ and any unit vector $\bm{u}$
we have 
\[\tc\tilde\gamma_{\bm{u}}\le\tc\gamma.\]
On the other hand total curvature of few developments gives an estimate 
for the total curvature of the original curve.
For example, if $\bm{i},\bm{j},\bm{k}$ is the orthonormal basis, then
\[\tc\gamma
\le
\tc\tilde\gamma_{\bm{i}}+\tc\tilde\gamma_{\bm{j}}+\tc\tilde\gamma_{\bm{k}}.\]

If $\gamma$ lies completely on the dark (or bright) side for direction $\bm{u}$,
then by Liberman's lemma we get 
\[\tc\tilde\gamma_{\bm{u}}\le \pi.\]
It follows that if $\gamma$ cross the horizons $\omega_{\bm{i}}$, $\omega_{\bm{j}}$ and $\omega_{\bm{k}}$
at most $N$ times, then 
\begin{align*}
\tc\gamma
&\le
\tc\tilde\gamma_{\bm{i}}
+\tc\tilde\gamma_{\bm{j}}
+\tc\tilde\gamma_{\bm{k}}
\le
\\
&\le(N+1)\cdot \pi.
\end{align*}
Therefore, 
to violate Main Theorem 
$\gamma$ has to cross the  horizons $\omega_{\bm{i}}$, $\omega_{\bm{j}}$ and $\omega_{\bm{k}}$ huge number of times.

\section{Curvature of development}\label{sec:curv-develop}

Let $\Sigma\subset\RR^3$
be a closed smooth strongly convex surface
and $\gamma\:[0,\ell]\to \Sigma$ be a unit-speed geodesic.
Assume that for some unit vector $\bm{u}$,
the geodesic $\gamma$ cross the horizon $\omega_{\bm{u}}$ transversely at 
$t_0<\dots <t_k$.
Set $\alpha_i=\measuredangle(\dot\gamma(t_i),\bm{u})-\tfrac\pi2$ for each $i$.
Note that $|\alpha_i|\le\tfrac\pi2$.

The values $t_i$ and  $\alpha_i$ 
will be called correspondingly \emph{meeting moments} 
and \emph{meeting angles}
of the geodesic $\gamma$ with the horizon $\omega_{\bm{u}}$.

Let us introduce new notation
\[\tc_{\bm{u}}\gamma\df\tc\tilde\gamma_{\bm{u}}.\]
From Liberman's lemma \ref{lem:liberman},
we get the following.

\begin{thm}{Corollary}\label{cor:liberman}
Let $\Sigma\subset\RR^3$
be a strongly convex smooth surface,
$\gamma\:[0,\ell]\to \Sigma$ be a unit-speed geodesic
and $\bm{u}$ is a unit vector.
Assume that  $\gamma$ cross the horizon $\omega_{\bm{u}}$ transversely
and 
$t_0<\dots <t_k$ be its meeting moments 
and $\alpha_0,\dots,\alpha_k$ be its meeting angles with the horizon $\omega_{\bm{u}}$.
Then
\[\tc_{\bm{u}}\gamma
\le 3\cdot\pi
+
2\cdot\left|\alpha_0-\alpha_1
+\dots +(-1)^k\cdot\alpha_k\right|.
\]

\end{thm}

As you will see further, 
in order to find the needed estimate the total curvature
of geodesic we will get an upper bound for the sum
\[\left|\alpha_0-\alpha_1
+\dots +(-1)^k\cdot\alpha_k\right|.\]
Finding  such an upper bound is the most important ingredient in the proof of the main theorem.

\parit{Proof.}
By Liberman's lemma,
\[\tc_{\bm{u}}(\gamma|_{[t_{i-1},t_i]})=\pm(\alpha_{i-1}-\alpha_i)\]
where the sign is ``$+$'' 
if  $\gamma_{[t_i,t_{i+1}]}$ lies on the dark side 
and ``$-$'' if it lies on the bright side from $\bm{u}$.
Summing all this up we get
\[\tc_{\bm{u}}(\gamma|_{[t_0,t_{k}]})
=
\left|\alpha_0
-2\cdot\alpha_1
+\dots+(-1)^{k-1}\cdot2\cdot\alpha_{k-1}+(-1)^k\cdot\alpha_k\right|.\]
By Liberman's lemma we also have
\[\tc_{\bm{u}}(\gamma|_{[0,t_0]}),\tc_{\bm{u}}(\gamma|_{[t_k,\ell]})\le \pi \]
Since $\alpha_0,\alpha_k\le\tfrac\pi2$, the statement follows.
\qeds 

If $\Sigma$ is a surface in $\RR^3$ and $p\in\Sigma$
we denote by $K_p$ the Gauss curvature of $\Sigma$ at $p$.

\begin{wrapfigure}{r}{37 mm}
\begin{lpic}[t(-5 mm),b(0 mm),r(0 mm),l(0 mm)]{pics/tongue(1)}
\lbl[t]{34,12;{$\omega_{\bm{u}}$}}
\lbl[lb]{26,24;{$\gamma$}}
\lbl[t]{1,11;{\tiny $\gamma(a)$}}
\lbl[t]{15,11;{\tiny $\gamma(b)$}}
\end{lpic}
\end{wrapfigure}
 
Assume $a,b$ be the meeting moments of minimizing geodesic $\gamma$ with $\omega_{\bm{u}}$.
The arc $\gamma|_{[a,b]}$ will be called \emph{$\omega_{\bm{u}}$-tongue}
if there is an immersed disc 
$\iota\:\DD\looparrowright\Sigma$ 
such that the closed curve $\iota|_{\partial D}$ is formed by joint of the arc $\gamma|_{[a,b]}$ and an arc of $\omega_{\bm{u}}$.
In this case the immersion $\iota$ is called the \emph{disc of the tongue}.

\begin{thm}{Tongue Lemma}\label{lem:tongue}
Let $\bm{u}$ be a unit vector,
$\gamma\:[a,b]\to\Sigma$ 
be a minimizing geodesic and $\omega_{\bm{u}}$-tongue on the strongly convex surface $\Sigma\subset\RR^3$ 
and 
$\iota\:\DD\looparrowright\Sigma$
be its disc.

Then 
\[\int\limits_{\DD} K_{\iota(x)}\cdot d_{\iota(x)}\area_\Sigma\] takes one of the value
\[\alpha-\beta,
\  -\alpha+\beta,
\ \pi-\alpha-\beta,
 \pi+\alpha+\beta \pmod{2\cdot\pi},
\]
where $\alpha$ and $\beta$ are the meeting angles at $a$ and $b$ correspondingly.

In particular 
\[\bigl|\alpha-\beta\bigr|
\le
\int\limits_{\DD} K_{\iota(x)}\cdot d_{\iota(x)}\area_\Sigma.
\eqlbl{eq:key2}\]
If in addition the image $\iota(\DD)$ lies  completely in the dark or bright side for $\bm{u}$, then
\[\tc_{\bm{u}}\gamma
\le 
\int\limits_{\iota(\DD)} K_p\cdot d_p\area_\Sigma<2\cdot \pi.\eqlbl{eq:key3}\]
\end{thm}

\parit{Proof.}
Since $\gamma$ is a geodesic, 
the parallel translation along $\gamma$ 
sends $\dot\gamma(a)$ to $\dot\gamma(b)$.

Note also that $\bm{u}$ belongs to the tangent plane to $\Sigma$ at any point on the horizon $\omega_{\bm{u}}$;
in particular the $\bm{u}$ extends to a parallel tangent vector field on $\omega_{\bm{u}}$.

It follows that parallel translation along $\iota|_{\partial\DD}$ 
rotates the tangent plane by angle 
\[\pm(\tfrac\pi2+\alpha)\pm(\tfrac\pi2+\beta)\]
To prove the main statement of the lemma,
it remains to apply Gauss--Bonnet formula.

Denote by $R$ the right hand side in \ref{eq:key2}.
Note that $R\ge 0$ 
and $|\alpha|,|\beta|\le \tfrac\pi2$.
From the main statement of the lemma it follows then that the minimal possible value for $R$ is $\bigl|\alpha-\beta\bigr|$.

To prove \ref{eq:key3}, note that 
in this case $\iota$ is an embedding.
Further note that the spherical image of the dark side of $\Sigma$ is hemisphere.
Therefore $2\cdot\pi$ is the integral of Gauss curvature along the dark side.
It follows that
\[\int\limits_{\DD} K_{\iota(x)}\cdot d_{\iota(x)}\area_\Sigma
=
\int\limits_{\iota(\DD)} K_p\cdot d_p\area_\Sigma
<2\cdot\pi.\]
By Liberman's lemma the statement follows.
\qeds

\section{Almost straight arcs}

Let $\eps>0$.
A curve $\gamma\:[a,b]\to\RR^3$ will be called \emph{$\eps$-straight}
if 
\[(1-\eps)\cdot\length \gamma\le |\gamma(b)-\gamma(a)|.\]

\begin{thm}{Lemma}\label{lem:eps-straight}
Assume $\eps>0$ 
and $n$ is a positive integer such that
$n\cdot\eps>2$.
Then any minimizing geodesic
on a convex surface $\Sigma$ in $\RR^3$
can be subdivided 
into $\eps$-straight arcs $\gamma_1,\dots,\gamma_n$.
\end{thm}

\parit{Proof.}
Let $\theta\in(0,\pi)$
be such that 
\[1-\cos\theta=\eps.\]

Assume two points $p$ and $q$ lie on the convex surface $\Sigma$.
Denote by $\nu_p$ and $\nu_q$ the outer normal vectors at $p$ and $q$ correspondingly.
Note that if 
\[\measuredangle(\nu_p,\nu_q)\le 2\cdot\theta,\]
then any minimizing geodesic from $p$ to $q$ on $\Sigma$
is $\eps$-straight.

Let $\gamma\:[0,\ell]\to \Sigma$ be a minimizing geodesic parametrized by its arc-length.

Assume $\gamma_{[t,\ell]}$ is not $\eps$-straight.
Set $t'$ to be the maximal value in $[t,\ell)$ such that the  interval $[t,t']$ is $\eps$-straight.
Consider a sequence $0=t_0<t_1<\dots<t_n<\ell$
such that $t_{i+1}=t_i'$ for each $i$.
Denote by $\nu_i$ the outer unit normal vector to $\Sigma$ at $\gamma(t_i)$. 
From above we get
\[\measuredangle(\nu_i,\nu_j)\ge2\cdot\theta\] for all $i\ne j$.
In other words, the open balls $\mathrm{B}_\theta(\nu_i)$ do not overlap 
in $\SS^2$.

It remains to note that 
\[\area[\mathrm{B}_\theta(\nu_i)] =2\cdot\pi\cdot\eps
 \ \ \text{and}\ \ 
 \area\SS^2=4\cdot\pi.
\]
Hence the result follows.
\qeds

\begin{thm}{Corollary}\label{cor:length-diameter}
Assume $\gamma\:[0,\ell]\to \Sigma$ is a unit-speed minimizing geodesic on the convex surface $\Sigma$ in $\RR^3$.
Then $\diam \gamma\ge \tfrac\ell{10}$.
\end{thm}

\parit{Proof.}
Apply Lemma~\ref{lem:eps-straight} for $\eps=\tfrac12$.
\qeds

\section{An arc in almost constant direction}

\begin{thm}{Proposition}\label{prop:almost-const}
For any $\eps>0$ there is $\delta>0$ such that the following holds.

If $\gamma\:[a,b]\to\Sigma$ is a minimizing geodesic 
on a smooth strongly convex surface $\Sigma$ in $\RR^3$,
then there is an interval $[a',b']\subset[a,b]$
such that 
\[\tc(\gamma|_{[a',b']})>\delta\cdot\tc\gamma.\]
and 
\[\measuredangle(\dot\gamma(t),\bm{u})<\eps\] 
for a fixed unit vector $\bm{u}$
and any $t\in[a',b']$.

Moreover, if $\eps=\tfrac1{10}$, then one can assume $\delta=\tfrac1{100^{100}}$.
\end{thm}

In the proof we will need the following two lemmas.

\begin{thm}{Lemma}\label{lem:almost-const-angles}
For any $\eps$ there is $\delta>0$ such that the following holds.

Assume $\gamma$ is a curve,
$\bm{v}_1$ and $\bm{v}_2$  are two vectors in $\RR^3$
and $0\le\alpha_1,\alpha_2\le\pi$ be such that
\begin{align*}
\eps
&<\measuredangle(\bm{v}_1,\bm{v}_2)<\pi-\eps
\\
\alpha_i-\delta
&<\measuredangle(\bm{v}_i,\dot\gamma(t))< \alpha_i+\delta.
\end{align*}
Then there is a vector $\bm{u}$ such that
$\measuredangle(\bm{u},\dot\gamma(t))<\eps$.

Moreover if $\eps<\tfrac1{10}$, 
then one can take $\delta=\eps^{10}$.
\end{thm}

The proof of the lemma above 
is straightforward computation;
we omit it.

\begin{thm}{Lemma}\label{lem:almost-const}
For any $\eps>0$ there is $\delta>0$
such that the following holds.

Let $\gamma\:[a,b]\to\Sigma$ 
be an $\delta$-straight minimizing geodesic 
on a smooth strongly convex surface $\Sigma$ in $\RR^3$.
Set $\bm{v}_\gamma=\gamma(b)-\gamma(a)$.
Then there in a sub-interval $[a',b']$ in $[a,b]$
such that 
\[\tc(\gamma|_{[a',b']})\ge\delta\cdot\tc\gamma.\]
and 
\[\alpha-\eps
\le
\measuredangle(\dot\gamma(t),\bm{v}_\gamma)
\le\alpha+\eps\] 
for some fixed $\alpha$
and
any $t\in[a',b']$.

Moreover if $\eps<\tfrac1{10}$ one can take $\delta=\eps^{10}$.
\end{thm}

\parit{Proof.}
Without loss of generality we can assume that
$a=0$, $b=2$ and 
\[\tc(\gamma|_{[1,2]})
\ge
\tfrac12\cdot\tc\gamma.\]

Set $p=\gamma(0)$.
Let $\theta\in(0,\pi)$ be such that $1-\cos\theta=\delta$.
Note that 
$$\measuredangle(\bm{v}_\gamma,\gamma(t)-p)
\le
\measuredangle(\tilde\gamma_p(1)-\tilde p,\tilde\gamma_p(2)-\tilde p)
<
\theta
\eqlbl{eq:angle1}$$
for any $t\ge 1$.

By Liberman's lemma 
\[\tc_p(\gamma|_{[1,2]})< \pi+\theta.\]
Assume $N=\lceil\tfrac\pi\theta+1\rceil$.
Then we can subdivide $\gamma|_{[1,2]}$ into $N$ arcs 
 $\gamma_1,\gamma_2,\dots,\gamma_N$ such that
\[\tc_p(\gamma_n)\le \theta\eqlbl{eq:tcp}\]
for each $n$.

From \ref{eq:angle1} and \ref{eq:tcp},
it follows that for each $n$, there is $\alpha_n$ such that
\[\alpha_n-\theta
\le
\measuredangle(\dot\gamma_n(t),\bm{v}_\gamma)
\le
\alpha_n+\theta.\] 
The arc $\gamma_n$ with the maximal total curvature will solve the proposition.

It remains to choose $\delta$ so that $\theta(\delta)<\tfrac\eps{100}$.
\qeds

\parit{Proof of Proposition~\ref{prop:almost-const}.}
Set $\gamma_0=\gamma$.

Fix $\delta>0$, set $n=\lceil\tfrac2\delta\rceil$.
By Lemma~\ref{lem:eps-straight}, the geodesic $\gamma_0$ can be subdivided into $n$ arcs which are  $\delta$-straight.
Let us choose the arc $\gamma'_0$ with the maximal total curvature.
Assuming $\delta<\tfrac1{10}$ we get
\[\tc\gamma'_0\ge\tfrac\delta{10}\cdot\tc\gamma_0.\]

Let $\alpha_1$ be the angle
and $\gamma_1$ be the arc in $\gamma'_0$ 
provided by Lemma~\ref{lem:almost-const}.
In particular 
\begin{align*}
\tc\gamma_1&\ge\delta\cdot\tc\gamma'_0
\ge
\\
&\ge \tfrac{\delta^2}{10}\cdot\tc\gamma_0.
\end{align*}

If 
$\alpha_1\le \tfrac\eps2$ or $\alpha_1\ge\pi-\tfrac\eps2$ 
and $\delta$ is sufficiently small,
then statement holds for the arc $\gamma_1$ and the vector $\bm{u}=\pm\bm{v}_{\gamma'_0}$.

Otherwise let us repeat the above construction for $\gamma_1$.
Namely, apply Lemma~\ref{lem:eps-straight} to the geodesic $\gamma_1$ and
denote by $\gamma_1'$ the $\delta$-straight 
arc with maximal total curvature.
If $\delta$ is small, we get 
\[
\tfrac\eps3
<
\measuredangle(\bm{v}_{\gamma_1'},\bm{v}_{\gamma'_0})
<
\pi-\tfrac\eps3
\eqlbl{eq:v1-perp-v_2}\]
Again, we get
\[\tc\gamma_1'\ge\tfrac\delta{10}\cdot\tc\gamma_1\ge \tfrac{\delta^3}{100}\cdot\tc\gamma_0.\]

Further apply Lemma~\ref{lem:almost-const} to $\gamma_1'$.
Denote by $\gamma_2$ and $\alpha_2$ the angle and the sub-arc of $\gamma_1'$.
Again
\[\tc\gamma_2\ge\tfrac{\delta^4}{100}\cdot\tc\gamma_0.\]

The curve $\gamma_2$ runs under nearly constant angle to  $\bm{v}_{\gamma'_0}$ and $\bm{v}_{\gamma'_1}$.
The inequality \ref{eq:v1-perp-v_2}
makes possible to apply Lemma \ref{lem:almost-const-angles}.
Hence the main statement in the proposition follows.

Straightforward computations prove the last statement.
\qeds

\section{Drifting geodesics}\label{sec:drifting}

In this section we fix notations which will be used further 
without additional explanation.

Fix a $(x,y,z)$-coordinates on the Euclidean space;
denote by $(\bm{i},\bm{j},\bm{k})$
the standard basis.

A plane parallel to say $(y,z)$-coordinate plane will be called $(y,z)$-plane.

\begin{thm}{Definition}
A smooth curve $\gamma\:[0,\ell]\to\RR^3$ 
is called \emph{$\bm{i}$-drifting} if both ends $\gamma(0)$ and $\gamma(\ell)$ lie on the $x$-axis and 
$\measuredangle(\dot\gamma(t),\bm{i})<\tfrac1{10}$ for all $t$.
\end{thm}

\parbf{$\bm{(\lambda,\mu,\nu)}$-frame.}
Let $\Sigma$ be a convex surface 
and $\gamma\:[0,\ell]\to \Sigma$ be an $\bm{i}$-drifting minimizing geodesic 
with unit-speed para\-me\-tri\-za\-tion.

Given $t\in [0,\ell]$, 
consider the oriented orthonormal frame $\lambda(t),\mu(t),\nu(t)$ 
such that $\nu(t)$ is the outer normal to $\Sigma$ at $\gamma(t)$,
the vector $\mu(t)$ is lies in $(y,z)$-plane 
and therefore the vector $\lambda(t)$ lies in the plane spanned by $\nu(t)$ and the $x$-axis.
We assume in addition that $\langle\lambda,\bm{i}\rangle\ge 0$.

Since $\langle\dot\gamma(t),\bm{i}\rangle>0$, w have
$\nu(t)\ne\bm{i}$.
Therefore the frame $(\lambda,\mu,\nu)$ is uniquely defined for any $t\in[0,\ell]$.

\parbf{Angle functions.}
Set 
\begin{align*}
\phi(t)&=\measuredangle(\bm{i},\dot\gamma(t)),&
\psi(t)&=\tfrac\pi2-\measuredangle(\bm{i},\nu(t)),&
\theta(t)&=\tfrac\pi2-\measuredangle(\mu(t),\dot\gamma(t)),
\end{align*}

From the above definitions it follows that $|\theta(t)|,|\psi(t)|\le \tfrac\pi2$ and for each $t$ there is a right spherical triangle with legs $|\theta(t)|,|\psi(t)|$ and hypotenuse $\phi(t)$.
In particular $\cos\theta\cdot\cos\psi=\cos\phi$.
Whence we get the following.

\begin{thm}{Claim}\label{clm:alpha-phi-psi}
For any $t$ we have 
\[
\phi(t)\ge |\psi(t)|\ \ \text{and}\ \ \phi(t)\ge |\theta(t)|
\]

\end{thm}

Applying Liberman's Lemma in the direction $\bm{i}$ we also get the following.

\begin{thm}{Claim}\label{clm:7.3}
If an arc $\gamma|_{[a,b]}$ lies in the dark (bright) side for $\bm{i}$, then 
the angle function $\phi$ is non-decreasing (correspondingly non-increasing)
in $[a,b]$.
\end{thm}

\section{Plane sections}

Assume $\gamma$ is a curve on a smooth strongly convex surface $\Sigma$ in $\RR^3$.
Consider a plane $L$ 
passing thru two points of $\gamma$, 
say $p=\gamma(a)$ and $q=\gamma(b)$ with $a<b$.
Let $L_\pm$ be a  half-planes in $L$ bounded by the line thru $p$ and $q$.
Set $\sigma_\pm=\Sigma\cap L_\pm$.

\begin{thm}{Observation}\label{obs:cut}
If $\gamma$ is a minimizing geodesic in the smooth strongly convex surface $\Sigma\subset \RR^3$ and $a$, $b$ and $\sigma_\pm$ as above, then
\[\length\sigma_\pm\ge\length(\gamma|_{[a,b]}).\]
\end{thm}

To prove the observation, it is sufficient to note that $\sigma_\pm$ are  smooth convex plane curve connecting $p$ to $q$ in $\Sigma$.

Based on this observation we give couple of estimates on drifting minimizing geodesics.

Let $\gamma\:[a,b]\to\RR^3$ be a curve and $\ell$ be a line which does not pass thru points of $\gamma$.
Assume $\phi\:[a,b]\to\RR$ be a continuous azimuth angle
of $\gamma$ in the cylindrical coordinates with the axis $\ell$.
If 
\[|\phi(b)-\phi(a)|\ge 2\cdot n\cdot\pi,\]
we will say that $\gamma$ \emph{goes around} the line $\ell$
at least $n$ times.

\begin{thm}{Proposition}\label{prop:around-once}
Assume $\gamma\:[0,\ell]\to \Sigma$ is an $\bm{i}$-drifting minimizing geodesic in the convex surface $\Sigma\subset \RR^3$, a sub-segment $[a,b]\subset [0,\ell]$ and the following conditions hold
\begin{enumerate}[(i)]
\item The points $\gamma(a)$ and $\gamma(b)$ lie in a half-plane with boundary line formed by the $x$-axis
and  the arc $\gamma|_{[a,b]}$ goes around the $x$-axis at least once.
\item  The $x$-coordinate of $\gamma(a)$ is larger that the $x$-coordinate of $\tfrac12\cdot(\gamma(0)+\gamma(\ell))$.
\end{enumerate}
Then  $\gamma(b)$ lies on the dark side for $\bm{i}$.
\end{thm}

\parit{Proof.} 
Let us apply Observation~\ref{obs:cut} to the plane containing $x$-axis, $\gamma(a)$ and $\gamma(b)$.

We can assume that $\gamma(0)$ is the origin of the $(x,y,z)$-coordinate system
and both points $p=\gamma(a)$ and $q=\gamma(b)$ lie in the $(x,z)$-coordinate half-plane with $x\ge 0$, denoted by $\Pi$.
We can assume that $\sigma_+\subset \Pi$.
Let $(x_p,0,z_p)$ and $(x_q,0,z_q)$ be the coordinates of $p$ and $q$.

From the assumptions we get $x_p<x_q<2\cdot x_p$.

Assume contrary,
then $\gamma(b)$ lies on the bright side for $\bm{i}$.
Then from convexity of the curve $\Pi\cap \Sigma$
we get 
\[\length\sigma_+\le \sqrt{(x_q-x_p)^2+z_p^2}.\]

On the other hand, since $\gamma|_{[a,b]}$ goes around $x$-axis at least once,
we get 
\[\length\gamma|_{[a,b]}\ge \sqrt{(x_q-x_p)^2+(z_p+z_q)^2}.\]

These two estimates contradict Observation~\ref{obs:cut}.
\qeds

\begin{thm}{Corollary}\label{cor:around-twice}
If $\Sigma$, $\gamma$, $\ell$, $a$ and $b$ as in the Proposition~\ref{prop:around-once} and the arc $\gamma|_{[a,b]}$ goes around the $x$-axis at least twice,
then the arc $\gamma|_{[b,\ell]}$ lies on the dark side with respect to $\bm{i}$.
\end{thm}

\parit{Proof.} 
Fix $b'\in [b,\ell]$.
Note that one can find $a'\in [a,b]$ 
such that the assumptions of Proposition~\ref{prop:around-once} hold for the interval $[a',b']$. 
Applying the proposition we get the result.
\qeds

\begin{thm}{Proposition}\label{prop:phi-psi}
Assume $\gamma\:[0,\ell]\to \Sigma$ is an $\bm{i}$-drifting minimizing geodesic 
in the convex surface $\Sigma\subset \RR^3$.
Assume that the arc $\gamma|_{[b,\ell]}$ 
lies on the dark side of $\Sigma$ with respect to $\bm{i}$.
If $b\le s<t\le \ell$ and the point $\gamma(s)$ lies in the plane $\Pi$ thru $\gamma(t)$ spanned by $\nu(t)$ and $\lambda(t)$,
then 
\[\phi(s)\le \psi(t).\]
\end{thm}

\parit{Proof.}
Let us apply Observation~\ref{obs:cut} to the plane $\Pi$ and $p=\gamma(s)$ and $q=\gamma(t)$.

Let $x_p$ and $x_q$ be the $x$-coordinates of $p$ and $q$.

Since $\gamma|_{[s,t]}$ lies in the dark side,
its Liberman's development $\tilde\gamma|_{[s,t]}$ 
with respect to $\bm{i}$ is concave.
In particular 
\[\length(\gamma|_{[s,t]})
=
\length(\tilde\gamma|_{[s,t]})
\ge
\tfrac{x_q-x_p}{\cos\phi(s)}.\]

On the other hand, convexity of $\sigma_+$ imply that
\[\length\sigma_+\le \tfrac{x_q-x_p}{\cos\psi(t)}.\]
It remains to apply Observation~\ref{obs:cut}.
\qeds

\section{\textit{s}-pairs}

Let $\Sigma\subset \RR^3$ 
be a strongly convex surface
and $\gamma\:[0,\ell]\to\Sigma$ be an $\bm{i}$-drifting minimizing geodesic.

Аssume that $\gamma$ intersects the horizon $\omega_{\bm{j}}$  transversely.

Let $t_0<t_1<\dots<t_k$ be the meeting moments of $\gamma$ with $\omega_{\bm{j}}$.
Set
\begin{align*}
\phi_n&=\phi(t_n)
&
\psi_n&=\psi(t_n)
&
\theta_n&=\theta(t_n).
\end{align*}
Note that  $\theta_n=\pm\alpha_n$
so we can say $s_n\cdot\theta_n=(-1)^n\cdot \alpha_n$ 
for some sequence of signs $s_i=\pm1$.
In particular
\[\alpha_0-\alpha_1+\dots+(-1)^k\cdot\alpha_k
=
s_0\cdot\theta_0+s_1\cdot\theta_1+\dots+s_k\cdot\theta_k.\]
Note that for the right choice of orientation,
if $s_n=+1$, 
then $\nu_{\gamma(t)}$ moves clockwise in $\SS^2$
at $t_n$
and if $s_n=-1$, then it moves counterclockwise.    

We say that a pair of indexes $i<j$
forms an \emph{$s$-pair} 
if 
\[
\sum_{n=i}^js_n=0\ \ 
\text{and}\ \ 
\sum_{n=i}^{j'}s_n>0
\]
if $i<j'<j$.

If you exchange ``$+1$'' and ``$-1$'' in $s$ by ``$($'' and ``$)$'' correspondingly, then $(i,j)$ is an $s$-pair
if and only if the $i$-th bracket forms a pair with $j$-bracket.

Note that any index $i$ appears in at most one $s$-pair and 
for any $s$-pair $(i,j)$ we have
\begin{itemize}
\item $s_i=1$; that is, $i$-th bracket has to be opening.
 \item $s_j=-1$; that is, $j$-th bracket has to be closing.
\end{itemize}
In particular,
\begin{align*}
s_i\cdot\theta_i+s_j\cdot\theta_j&=\theta_i-\theta_j=
\\
&=(-1)^i\cdot\alpha_i+(-1)^j\cdot\alpha_j.
\end{align*}

\begin{wrapfigure}{r}{52 mm}
\begin{lpic}[t(-7 mm),b(1 mm),r(0 mm),l(0 mm)]{pics/s-pair(1)}
\lbl[br]{2,4;$+$}
\lbl[br]{3,14;$+$}
\lbl[br]{4,24;$+$}
\lbl[br]{6,34;$+$}
\lbl[bl]{16,34;$-$}
\lbl[br]{28,34;$+$}
\lbl[br]{33,44;$+$}
\lbl[bl]{43,44;$-$}
\lbl[bl]{45,34;$-$}
\lbl[bl]{46,24;$-$}
\lbl[bl]{47,14;$-$}
\lbl[bl]{47.5,4;$-$}
\lbl[lt]{4,2;$\gamma(t_i)$}
\lbl[rt]{46,2;$\gamma(t_j)$}
\end{lpic}
\end{wrapfigure}

\parbf{Tongue interpretation.}
Assume $(i,j)$ is an $s$-pair.
Note that in this case there is an arc of $\omega_{\bm{j}}$
from $\gamma(t_i)$ to $\gamma(t_j)$
with monotonic $x$-coordinate.
Moreover a disc of the tone has this arc in the boundary.

The proof can be guessed from the diagram.
It shows a lift of $\gamma$ in the universal cover of strip of $\Sigma$ between $(y,z)$-planes containing $\gamma|_{[t_i,t_j]}$;
the solid horizontal lines correspond are lifts of $\omega_{\bm{j}}$.

We say that $q$ is the depth of an $s$-pair $(i,j)$
(briefly $q=\depth_{s}(i,j)$) 
if $q$ is the maximal number such that there is $q$-long nested sequence of $s$-pairs starting with $(i,j)$.
For example the $s$-pair on the diagram has depth $5$.

More precisely, the depth of $(i,j)$ is the maximal number $q$
for which there is a sequence of $s$-pairs
$(i,j)=(i_1,j_1),(i_2,j_2),\zz\dots,(i_q,j_q)$ such that
\[i=i_1<\dots<i_q<j_q<\dots<j_1=j.\]

Note that the $s$-pairs of the same depth do not overlap;
that is if  $\depth(i,j)=\depth(i',j')$
for two distinct $s$-pairs $(i,j)$ and $(i',j')$,
then either $i<j<i'<j'$ or $i'<j'<i<j$.

The following proposition follow directly from the discussion above.

\begin{thm}{Proposition}\label{prop:immersion}
Let $(i,j)$ be an $s$-pair.
Then the arcs $\gamma|_{[t_i,t_j]}$ and an arc of $\omega_{\bm{j}}$ bound an immersed disc in $\Sigma$ which lies between $(y,z)$-planes thru $\gamma(t_i)$ and $\gamma(t_j)$.
Moreover the maximal multiplicity of the disc is at most $\depth_{s}(i,j)$.
\end{thm}

\begin{thm}{Corollary}\label{cor:Sq}
Let us denote by $S_q$  the subset of indexes $\{1,\dots,k\}$
which are the parts of $s$-pairs with depth $q$.
Then
\[\left|\sum_{n\in S_q}(-1)^n\cdot\alpha_n\right|
=\left|\sum_{n\in S_q}s_n\cdot\theta_n\right|
\le 4\cdot\pi\cdot q.\]
\end{thm}

\parit{Proof.} For each $n$ denote by $K_n$ the integral of Gauss curvature of the part of surface $\Sigma$ with the $x$-coordinate less than the $x$-coordinate of $\gamma(t_n)$.
Note that 
\[0\le K_1\le\dots\le K_k\le 4\cdot\pi.\]

By Proposition~\ref{prop:immersion} and the Tongue Lemma,
we get
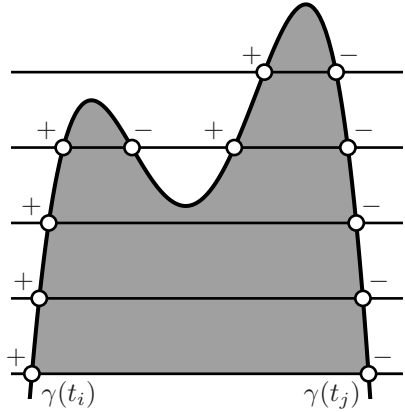
\[s_i\cdot\theta_i+s_j\cdot\theta_j=\theta_i-\theta_j\le q\cdot (K_j-K_i)\]

The statement follows since the $s$-pairs with the same depth do not overlap.
\qeds

\begin{thm}{Corollary}\label{cor:gamma-0}
Assume 
\[q=\max_{1\le i<j\le k}\left\{\biggl|\sum_{n=i}^js_n\biggr|\right\}\]
Then
\[\left|\sum_{n=1}^k s_n\cdot\theta_n\right|
\le 2\cdot q\cdot(q+\tfrac32)\cdot \pi.
\]
\end{thm}

\parit{Proof.}
Denote by $S$ the set of all indexes which appear in some $s$-pair.

Note that depth of any $s$-pair is at most $q$.
That is,
\[S=S_1\cup\dots\cup S_q.\]
By Corollary~\ref{cor:Sq},
\[\left|\sum_{n\in S}s_n\cdot\theta_n\right|\le 2\cdot q\cdot(q+1)\cdot\pi.\eqlbl{sum-S}\]

Set $R=\{1,\dots,k\}\backslash S$;
this is the set of indexes 
which do not appear in an $s$-pair.

Given $r$, set $i\in Q_r$ 
if
\[\sum_{n=1}^is_n=r.\]
Note that $Q_r\ne\emptyset$ for at most $q$ values of $r$
and in each set $Q_r$ there are at most $2$ indexes 
which do not appear in an $s$-pair;
that is $Q_r\cap R$ has at most two indexes for each $r$.

Since $|a_n|<\tfrac\pi2$, we get
\[\left|\sum_{n\in R}s_n\cdot\theta_n\right|
\le
q\cdot\pi.
\]
The latter inequality together with \ref{sum-S} implies the statement in the corollary.
\qeds

\section{Geometric growth}\label{sec:geometric-growth}

\begin{thm}{Claim}\label{clm:alpha-psi}
Assume $\psi(t)>\eps$ for  $t\in[t_{i},t_{i+1}]$
and $s_i=s_{i+1}$
Then 
\[|\theta_{i+1}-\theta_i|>\pi\cdot\sin\eps.\]

\end{thm}

\parit{Proof.}
Note that the arc $\gamma|_{[t_{i},t_{i+1}]}$ is a tongue with embedded disc $\iota\:\DD^2\to \Sigma$.
Since $\psi(t)>\eps$, 
the spherical image $\nu\circ\iota(\DD^2)$ 
of $\iota(\DD^2)$ 
lies in a half-disc of radius $\tfrac\pi2-\eps$ in $\SS^2$.
Note that
\[K(\iota(\DD^2))
=
\area(\nu\circ\iota(\DD^2))
<
\pi\cdot(1-\sin\eps).\]
It remains to apply Tongue Lemma~\ref{lem:tongue}.
\qeds

\begin{thm}{Claim}\label{clm:geometric-grouth}
Assume $\gamma$ lies on the dark side for $\bm{i}$.
Then for any pair of indexes $j>i$,
such that 
\[\biggl|\sum_{n=i}^{j}s_n\biggr|> 5\]
we have
\[\phi_j>\tfrac32\cdot\phi_i.\]

\end{thm}

\parit{Proof.}
By Claim~\ref{clm:7.3}, we may assume that 
\[\sum_{n=i}^{j}s_n= 6\]

Let $j'$ be the least index 
such that
\[\biggl|\sum_{n=i}^{j'} s_n\biggr|=5.\]

Note that for any $b>t_j$ there is $a\in[t_i,t_j]$
such that interval $[a,b]$ satisfies the assumptions of Proposition~\ref{prop:phi-psi}.
In particular $\psi(b)>\phi_i$ for any $b>t_j$.
Applying Claim~\ref{clm:alpha-psi},
we get that $|\theta_j|>\tfrac\pi2\cdot \phi_i$ or 
$|\theta_{j'}|>\tfrac\pi2\cdot \phi_i$.
By Claim~\ref{clm:7.3}, $\phi_n$ is non-decreasing,
and $\phi_n\ge |\theta_n|$ for any $n$,
in both cases we get
\[\phi_j>\tfrac\pi2\cdot \phi_i.\]
Hence the result follows.
\qeds

\begin{thm}{Proposition}\label{prop:graph}
If $\gamma$ is an $\bm{i}$-drifting minimizing geodesic on the dark side for $\bm{i}$, then
\[\tc_{\bm{j}}\gamma\le 100\cdot\pi.\]

\end{thm}

\parit{Proof.}
We can assume that
$\gamma$ cross the $\bm{j}$-horizon $\omega_{\bm{j}}$ transversely.
Let $t_0<\dots<t_k$ be the meeting moments of $\gamma$ with $\omega_{\bm{j}}$ and
$s_0,\dots,s_k$ the signs.

Recall that $S_q$ denotes the subset of indexes $\{1,\dots,k\}$
which appear in $s$-pair with depth $q$.
According to Corollary~\ref{cor:Sq},
\[\left|\sum_{n\in S_q}s_n\cdot\theta_n\right|\le 4\cdot q\cdot \pi.\]
In particular,
\[\left|\sum_{n\in S_1\cup\dots\cup S_5}
s_n\cdot\theta_n\right|
\le 
40\cdot\pi.\]

Set $R=\{1,\dots,k\}\backslash (S_1\cup\dots\cup S_5)$;
this is the set of indexes which appear in $s$-pairs with depth at least $6$ 
as well as those which do not appear in any $s$-pair.

According to Claim~\ref{clm:alpha-phi-psi},
\[\left|\sum_{n\in R}
s_n\cdot\theta_n\right|
\le
\sum_{n\in R}
|\theta_n|\le \sum_{n\in R}\phi_n.\] 
To estimate the last sum will use the results in Section~\ref{sec:geometric-growth}.
First let us subdivide $R$ into 5 subsets $R_1,\dots,R_5$,
by setting 
$n\in R_m$ if $m\equiv n\pmod 5$.

Given $n\in R_m$, denote by $n'$ the least index in $R_m$ which is larger $n$;
$n'$ is defined for any $n\in R_m$ except the largest one.
According to Claim~\ref{clm:geometric-grouth}, 
$\phi_{n'}>\tfrac32\cdot \phi_n$;
that is, the sequence $(\phi_n)_{n\in R_m}$ grows faster than the geometric progression with coefficient $\tfrac32$.
Since $\phi_n$ is non-decreasing in $n$,
we get 
\[\sum_{n\in R_m}\phi_n< 3\cdot\phi_k.\]
It follows that 
\[\sum_{n\in R}\phi_n< 15\cdot\phi_k\le\tfrac{15}2\cdot\pi.\]

By Corollary~\ref{cor:liberman},
\begin{align*}
\tc_{\bm{j}}\gamma
&\le 
2\cdot\pi
+
2\cdot [\alpha_0-\alpha_1+\dots+(-1)^k\cdot\alpha_k]< 
\\
&< 100\cdot \pi.
\end{align*}
\qeds

\section{Assembling of the proof}

Assume $\gamma\:[0,\ell]\to \Sigma$ is a minimizing geodesic in a convex surface $\Sigma\subset \RR^3$.

According to propositions \ref{prop:semicontinuity} and \ref{prop:convegence} we can assume that $\Sigma$ is closed,  strongly convex and smooth
and the geodesic $\gamma$ has finite length.

According to Proposition~\ref{prop:almost-const}, 
we can pass to an $\bm{i}$-drifting arc $\gamma'$
of $\gamma$ for some $(x,y,z)$-coordinate system 
such that 
\[\tc\gamma'>\tfrac{1}{100^{100}}\cdot\tc\gamma.
\eqlbl{eq:1/100}
\]
We will use the notations in Section~\ref{sec:drifting} for $\gamma'$.

Rotating $(y,z)$-coordinate plane we can ensure that
\[\tc\gamma'\le 10\cdot\tc_{\bm{j}}\gamma'\]
and that $\gamma'$ cross the horizon $\omega_{\bm{j}}$ transversally.

By Corollary \ref{cor:around-twice},
we can subdivide $\gamma'$ into at most three arcs: 
\begin{itemize}
\item \emph{Left arc} $\gamma_-'$ which lies on the bright side for $\bm{i}$,
\item \emph{Middle arc} $\gamma_0'$ which rotates around $x$-axis at most $4$ times.
\item \emph{Right arc} $\gamma_+'$ which lies on the dark side for $\bm{i}$.
\end{itemize}

Indeed, choose an arc $\gamma'|_{[a,b]}$  
on the right from the $(y,z)$-plane thru
$\tfrac12\cdot(\gamma'(0)+\gamma'(\ell))$
which rotates around $x$-axis 2 times 
and assume that $b$ takes the minimal possible value.  
Note that if $\gamma'(s)$ lies on $(y,z)$-plane thru $\tfrac12\cdot(\gamma'(0)+\gamma'(\ell))$,
then $[s,b]\supset[a,b]$ 
and any sub-arc of $[s,b]$ rotates around $x$-axis at most 2 times.

By Corollary \ref{cor:around-twice}, 
we can take $\gamma'_+=\gamma'|_{[b,\ell]}$;
in case if there is no such arc $[a,b]$, we assume that $\gamma_+'$ is not presented.
Repeat the construction reverting the direction of $x$-axis;
we get the left arc $\gamma_-'$.
The remaining arc is assumed to be $\gamma_0'$; note that any sub-arc of $\gamma_0'$ is divided by the $(y,z)$-plane thru $\tfrac12\cdot(\gamma'(0)+\gamma'(\ell))$ into two each of which rotate around $x$-axis at most 2 times. 
Therefore the number of rotations of any arc in $\gamma_0'$ is at most $4$.

Let us estimate the total curvature of $\gamma'_-$, $\gamma'_0$ and $\gamma'_+$ separately.

By Proposition~\ref{prop:graph}, we get 
\[\tc_{\bm{j}}\gamma'_+\le 100\cdot\pi.\eqlbl{+}\]
Similarly  
\[\tc_{\bm{j}}\gamma'_-\le 100\cdot\pi.\eqlbl{-}\]

By Corollary~\ref{cor:gamma-0},
\[\tc_{\bm{j}}\gamma'_0\le 100\cdot\pi.\eqlbl{0}\]
Together with \ref{+}, \ref{-} and \ref{0}
the latter implies that 
\[\tc_{\bm{j}}\gamma'\le 300\cdot\pi.\]
From \ref{eq:1/100} the result follows.
\qeds

\section{Final remarks}

Note that our main theorem does not hold if one removes word \emph{minimizing} from its formulation, even if the length of geodesic is bounded by the intrinsic diameter of the surface.

In addition to Corollary~\ref{cor:length-diameter}, 
we know couple of properties of the minimizing geodesics which distinguish them in the class of all geodesics on a convex surfaces.
We want to mention them altho they did not help us in the proof.
The first of these properties was essentially discovered by Anatolii Milka in \cite{milka-bending}.
To keep things simpler, we give a formulation only in the smooth case.

\begin{thm}{Stretching lemma}
Let $\gamma\:[0,\ell]\to\Sigma$ be a unit-speed minimizing geodesic from $p$ to $q$ on smooth convex hypersurface $\Sigma\subset\RR^m$ and $p_t=\gamma(t)-\dot\gamma(t)\cdot t$.
Then for any $t\in [0,\ell]$, the point $q$ lies on the dark side from $p_t$.
\end{thm}

\begin{thm}{Convex hat lemma}
Let $\Sigma\subset\RR^m$ be a closed convex hypersurface and
$\Pi$ be a hyperplane which cuts $\Sigma$ into two parts $\Sigma_1$ and $\Sigma_2$.
Assume that the reflection of $\Sigma_1$ in $\Pi$ lies inside $\Sigma$.
Then $\Sigma_1$ forms a convex set in $\Sigma$;
that is, any minimizing geodesic in $\Sigma$ with the ends in $\Sigma_1$ lies completely in $\Sigma_1$.
\end{thm}

The proofs of both statements can be found in \cite{petrunin-orthodox}.

\medskip

It is ashamed to confess, but we were not able to generalize the main theorem to higher dimensions.
Namely we can not give an answer to the following question.

\begin{thm}{Open question}
Is it true that for any positive integer $m$,
there is a constant $C_m$
such that the total curvature of arbitrary minimizing geodesic on a convex hypersurface in $\RR^m$ does not exceed $C_m$?
\end{thm}

\Addresses

\begin{thebibliography}{99}
\bibitem{AH-PSV}
Agarwal, P. K.; Har-Peled, S.; Sharir, M.; Varadarajan, K. R.,
``Approximating shortest paths on a convex polytope in three dimensions''.
\textit{J. ACM}
44.4 (1997),
pp. 567--584.
\bibitem{pach}
Pach, J.,
``Folding and turning along geodesics in a convex surface'',
\textit{Geombinatorics}
7.2 (1997)
pp. 61--65.
\bibitem{BKZ}
B{\'a}r{\'a}ny, I.; Kuperberg, K.; Zamfirescu, T.,
``Total curvature and spiralling shortest paths''.
\textit{Discrete Comput. Geom.}
30.2 (2003),
pp. 167--176.
\bibitem{liberman}
Либерман И. М.,
«Геодезические
 линии
 на
 выпуклых
 поверхностях»,
\textit{ДАН СССР}
32.5 (1941),
с. 310---313. 

\bibitem{usov}
Усов, В. В. 
«О длине сферического изображения геодезической на выпуклой поверхности». \textit{Сибирский математический журнал} 
17.1 (1976), 
с. 233---236.

\bibitem{berg}
Berg, I. D. 
"An estimate on the total curvature of a geodesic in Euclidean 3-space-with-boundary." 
\textit{Geometriae Dedicata} 
13.1 (1982),
pp. 1--6.

\bibitem{pogorelov}
Погорелов, А. В., 
\textit{Внешняя геометрия выпуклых поверхностей.} 
1969.

\bibitem{zalgaller}
Залгаллер, В. А. 
«Вопрос о сферическом изображении кратчайшей».
\textit{Укр. геометрический сб.}
10 (1971) 
с. 12---18.

\bibitem{milka}
Милка, А. Д. 
«Кратчайшая с неспрямляемым сферическим изображением». 
\textit{Укр. геометрический сб.} 16 (1974)
с. 35---52.

\bibitem{usov-conj-pog} 
Усов, В. В. 
«О пространственном повороте кривых на выпуклых поверхностях». 
\textit{Сибирский математический журнал}
17.6 (1976),
с. 1427---1430.

\bibitem{milka-liberman}
Милка, А. Д. 
«Аналог теоремы Либермана в римановом пространстве».
Украинский геометрический сборник,
24 (1981), 
с. 82---84,

\bibitem{petrunin}
Petrunin, A.
``Applications of quasigeodesics and gradient curves''.
\textit{Math. Sci. Res. Inst. Publ.}
30 (1997),
pp. 203--219

\bibitem{milka-bending}  Милка, А. Д. «Кратчайшие  линии  на  выпуклых  поверхностях».  Докл.  АН  СССР,   1979,   248, \textnumero1,  34---36  

\bibitem{petrunin-orthodox}
Petrunin, A.
``Puzzles in geometry that I know and love''.
\texttt{arXiv:0906.0290 [math.HO]}
\end{thebibliography}
\end{document}


\title{О вариации поворота кратчайших\\
 на выпуклых поверхностях}
\author{Нина Лебедева и Антон Петрунин}
\nofootnote{Работа
выполнена при финансовой поддержке первого из авторов фондом РФФИ грант 14-01-00062, второго фондом NSF, грант DMS 1309340.}

\newcommand{\Addresses}{{\bigskip\footnotesize

  Нина Лебедева, \par\nopagebreak\textsc{198504, Санкт-Петербург, Старый Петергоф, Университетский пр., дом 28.}
  \par\nopagebreak
  \textsc{191023, Санкт-Петербург, наб. Фонтанки, 27.}
  \par\nopagebreak
  \textit{Электронная почта}: \texttt{lebed@pdmi.ras.ru}

\medskip

  Anton Petrunin, \par\nopagebreak\textsc{Math. Dept. PSU, University Park, PA 16802, USA}
  \par\nopagebreak
  \textit{Электронная почта}: \texttt{petrunin@math.psu.edu}
}}
\date{}
\maketitle

\begin{abstract}
Показано, что вариация поворота кратчайших
на выпуклой поверхности в трёхмерном евклидовом пространстве 
не превосходит некоторой константы.
\end{abstract}

\section{Введение}

Эта заметка отвечает на вопрос заданный нам Дмитрием Бураго и
этот же вопрос мы нашли позднее в \cite{AH-PSV}, \cite{pach} и \cite{BKZ}. 

\begin{thm}{Основная теорема}\label{thm:main}
Вариация поворота кратчайшей
на выпукой поверхности в трёхмерном евклидовом пространстве $\RR^3$ 
не превосходит $1000^{1000}$.
\end{thm}

Значение $2\cdot\pi$ даёт точную оценку для аналогичной задачи на плоскости.
Есть примеры кратчайших на выпукой поверхности в $\RR^3$
с вариацией поворота строго больше $2\cdot\pi$ и
предположительно, точная оценка только немного превосходит $2\cdot\pi$.
Такой пример был построен Барань, Куперберг и Замфиреску в \cite{BKZ}. 
Вот ещё несколько связанных результатов.

\begin{itemize}
\item В \cite{liberman}, Либерман даёт оценку на вариацию поворота короткой геодезической на поверхности компактного тела $K$ через отношение диаметра $K$ и радиуса вписанного в $K$ шара.
Его доказательство использует аналог леммы \ref{lem:liberman} ниже.

\item В \cite{usov}, 
Усов даёт точную оценку на вариацию поворота геодезической на графике $\ell$-липшицевой выпуклой функции. 
А именно, он доказывает, что если $f\:\RR^2\to\RR$ является $\ell$-липшицевой и выпуклой тогда любая геодезическая на графике \[\Gamma_f=\set{(x,y,z)\in \RR^3}{z=f(x,y)}\] 
имеет вариацию поворота не больше $2\cdot \ell$.
Забавное обобщение этого результата дано Бергом в \cite{berg}.

\item В \cite{pogorelov}, Погорелов сформулировал гипотезу о том, что сферический образ геодезической на выпуклой поверхности обязан быть (локально) спрямляемым.
Легко видеть, что длина сферического образа геодезической не может быть меньше её вариации поворота.
Таким образом, эта гипотеза может рассматриваться как усиление теоремы Либермана сформуливанной выше.
Контрпримеры к различным вариантам этой гипотезы были получены Залгаллером в \cite{zalgaller},
Милкой в \cite{milka}
и Усовым в \cite{usov-conj-pog}.
Часть этих результатов была переоткрыта Пахом в \cite{pach}.

\item В \cite{BKZ},
Барань,
Куперберг
и Замфиреску построили винтовую кратчайшую на выпуклой поверхности;
то есть кратчайшую, которая обходит вокруг данной прямой произвольно большое число раз.
Ими были также переоткрыты результаты Либермана и Усова упомянутые выше.
\end{itemize}

\parbf{Идея доказательства.}
Задача сводится к случаю, когда кратчайшая $\gamma$ идёт приблизительно в одном направлении, 
скажем $\gamma$ параметризована длиной и $\dot\gamma(t)\approx\bm{i}$.

В этом случае, для оценки вариации поворота достаточно оценить интеграл
\[\int\langle\ddot\gamma(t),\bm{j}\rangle\cdot dt\] 
для вектора $\bm{j}\perp\bm{i}$.
Представьте, что поверхность освещена в направлении вектора $\bm{j}$ и, таким образом, некоторая кривая $\omega$ разделяет её на тёмную сторону и светлую.
На рисунке изображены несколько комбинаторных вариантов, как $\gamma$ может встречаться с $\omega$.

\begin{center}
 \begin{lpic}[t(2 mm),b(2 mm),r(0 mm),l(0 mm)]{pics/spiral-and-snake(1)}
\lbl[b]{93,14;$\omega$}
\end{lpic}
\end{center}

В первом случае вариация поворота оценивается через интеграл Гауссовой кривизны регионов заключённых между $\gamma$ и $\omega$.
Эта оценка следует из \emph{леммы о языке}~\ref{lem:tongue},
которая является основным моментом доказательства.

Второй случай может показаться невозможным, 
однако именно так ведёт себя винтовая кратчайшая построенная в  \cite{BKZ}.
В этом случае мы считаем вариацию поворота для каждого оборота кратчайшей и доказываем, что полученная последовательность чисел растёт быстрее чем геометрически от некоторого места в середине к обоим концам.
Оценка следует поскольку вариация поворота последних оборотов не превосходит $2\cdot\pi$.

Последний случай может рассматриваться как смесь первых двух
и его доказательство проводится смешанной техникой.

\section{Предварительные сведения}

\parbf{Вариация поворота.}
\emph{Вариация поворота}, $\tc\gamma$, кривой $\gamma\:[0,\ell]\to \RR^3$ 
определяется как точная верхняя грань сумм внешних углов
ломанных вписанных в $\gamma$.

Заметим, что вариация поворота ломанной $\sigma$ равна сумме её внешних углов.

Если $\gamma$ гладкая кривая параметризованная длиной, 
то 
\[\tc\gamma=\int\limits_0^\ell \kappa(t)\cdot dt,\]
гдя $\kappa(t)=|\ddot\gamma(t)|$ обозначает кривизну $\gamma$ в момент $t$.

\begin{thm}{Предложение}\label{prop:semicontinuity}
Предположим $\gamma_n\:\II\to\RR^3$ есть последовательность кривых, сходящаяся поточечно к кривой $\gamma_\infty\:\II\to\RR^3$.
Тогда 
\[\liminf_{n\to\infty}\tc\gamma_n\ge \tc\gamma_\infty.\]
\end{thm}

\parit{Доказательство.}
Выберем ломанную $\sigma_\infty$ вписанную в $\gamma_\infty$.
Пусть $\gamma_\infty(t_0),\zz\dots,\gamma_\infty(t_k)$
её вершины в порядке появления на $\gamma_\infty$.
Рассмотрим ломанную $\sigma_n$ вписанную в $\gamma_n$ с вершинами 
$\gamma_n(t_0),\zz\dots,\gamma_n(t_k)$.
Заметим, что
\[\tc\sigma_n\to\tc\sigma_\infty\ \ \text{при}\ \ n\to\infty.\]
По определению вариации поворота, 
\[\tc\sigma_n\le\tc\gamma_n.\]
Предложение следует, поскольку ломаная $\sigma_\infty$ 
может быть выбрана так, что 
$\tc\sigma_\infty$ произвольно близко к $\tc\gamma_\infty$.
\qeds

\parbf{Сходимость множеств.}
Для замкнутого множества $\Sigma\subset \RR^3$,
обозначим через $\dist_\Sigma$ функцию расстояния до $\Sigma$,
то есть
\[\dist_\Sigma x=\inf\set{|x-y|}{y\in\Sigma}.\]

Последовательность $\Sigma_n\subset \RR^3$ замкнутых множеств 
\emph{сходится} к замкнутому множеству $\Sigma_\infty\subset \RR^3$ 
если 
$\dist_{\Sigma_n}x\to\dist_{\Sigma_\infty}x$ для любого $x\in\RR^3$ при $n\to\infty$.

\parbf{Выпуклые поверхности.}
\emph{Выпуклая поверхность} $\Sigma$ в трёхмерном евклидовом пространстве  $\RR^3$ определяется как граница замкнутого выпуклого множества $K$ с непустой внутренностью.
Если $K$ компактно, мы говорим, что $\Sigma$ \emph{замкнутая}.

Предположим $\Sigma$ гладкая.
Если в любой точке $\Sigma$ 
её главные кривизны положительны, 
то $\Sigma$ называется \emph{сильно выпуклой}. 

\begin{thm}{Предложение}\label{prop:convegence}
Пусть $\Sigma_n$ есть последовательность выпуклых поверхностей, сходящихся к выпуклой поверхности $\Sigma_\infty$.
Тогда для любой кратчайшей $\gamma_\infty$ на $\Sigma_\infty$ 
существует последовательность кратчайших $\gamma_n$ на $\Sigma_n$, сходящаяся поточечно к $\gamma_\infty$ при $n\to \infty$.
\end{thm}

\parit{Доказательство.}
Предположим кратчайшая $\gamma_\infty\:[0,\ell]\to\Sigma_\infty$ параметризована длиной.

Выберем интервал $[a,b]\subset (0,\ell)$.
Пусть $p_\infty=\gamma_\infty(a)$ и $q_\infty=\gamma_\infty(b)$.
Выберем две последовательности точек $p_n,q_n\zz\in\Sigma_n$, сходящихся к $p_\infty$ и $q_\infty$ соответственно.

Обозначим через $\gamma_n$ кратчайшую между $p_n$ и $q_n$ на $\Sigma_n$.
Заметим, что $\gamma_n$ сходится к дуге $\gamma_\infty|_{[a,b]}$
при $n\to\infty$.

Выбирая интервал $[a,b]$ ближе и ближе к $[0,\ell]$ и переходя к диагональной подпоследовательности кратчайших получаем требуемое утверждение.
\qeds

\section{Лемма Либермана}

В этом разделе мы дадим небольшое обобщение конструкции Либермана из \cite{liberman};
см. также \cite{milka-liberman} и \cite{petrunin}.

\parbf{Развёртка.}
Пусть $\gamma\:[0,\ell]\to\RR^3$ есть кривая параметризованная длиной и точка $p$ не лежит на $\gamma$.

Предположим $\tilde\gamma_p\:[0,\ell]\to\RR^2$ есть плоская кривая параметризованная длиной и для 
некоторой точки $\tilde p$ на плоскости выполняется   
\[|\tilde p-\tilde\gamma(t)|=|p-\gamma(t)|\]
для всех $t\in[0,\ell]$.
Более того, направление из $\tilde p$ в $\tilde \gamma(t)$ монотонно крутится (по часовой или против часовой стрелке). 
Тогда $\tilde\gamma_p$ называется \emph{развёрткой} $\gamma$ относительно $p$.

\begin{center}
\begin{lpic}[t(0 mm),b(9 mm),r(0 mm),l(0 mm)]{pics/convex-concave(1)}
\lbl[tl]{23,16;$\tilde\gamma_p(a)$}
\lbl[tr]{4,16;$\tilde\gamma_p(b)$}
\lbl[tl]{14,1;$\tilde p$}
\lbl[ll]{67,18;$\tilde\gamma_p(a)$}
\lbl[rt]{45,18;$\tilde\gamma_p(b)$}
\lbl[tl]{57,1;$\tilde p$}
\lbl[t]{13,-5;Выпуклая развёртка.}
\lbl[t]{56,-5;Вогнутая развёртка.}
\end{lpic}
\end{center}

Развёртка $\tilde\gamma_p$ называется \emph{выпуклой} (\emph{вогнутой}) на интервале $[a,b]$
если дуга $\tilde\gamma_p|_{[a,b]}$
отрезает от угла $\angle \tilde p^{\tilde\gamma_p(a)}_{\tilde\gamma_p(b)}$
выпуклую ограниченную (соответственно неограниченную) область.

Мы говорим, что $\tilde\gamma_p$ \emph{локально} выпуклая (вогнутая) на интервале $[a,b]$
если любая точка $x\in [a,b]$ допускает замкнутую окрестность $[a',b']$ в $[a,b]$
такую, что $\tilde\gamma_p$ выпукла (соответственно вогнута) на интервале $[a',b']$.

Перейдя к пределу приведённой конструкции при $p$ убегающем вдоль луча в направлении единичного вектора $(-\bm{u})$, мы получим предельную плоскую кривую называемую развёрткой $\gamma$ в направлении  $\bm{u}$, она будет обозначаться $\tilde\gamma_{\bm{u}}$.

Развёртка $\tilde\gamma_{\bm{u}}$ также определяется следующими свойствами
(1) развёртка $\tilde\gamma_{\bm{u}}\:[0,\ell]\to \R^2$ параметризована длиной,
(2) для некого единичного вектора $\tilde{\bm{u}}\in\RR^2$ равенство
\[\langle \tilde{\bm{u}},\tilde\gamma_{\bm{u}}(t)\rangle
=
\langle  \bm{u},\gamma(t)\rangle\]
выполняется при всех $t\in [0,\ell]$ 
и
(3) проекция $\tilde\gamma_{\bm{u}}(t)$ на прямую перпендикулярную к $\tilde{\bm{u}}$
монотонна по $t$.

Без ущерба для общности, можно предположить, что $\tilde{\bm{u}}$ является вертикальным вектором координатной плоскости.
В этом случае $\tilde\gamma_{\bm{u}}$ выпукла (вогнута) на интервале $[a,b]$ 
если двуугольник ограниченный дугой $\tilde\gamma_{\bm{u}}|_{[a,b]}$
и отрезком $[\tilde\gamma_{\bm{u}}(a)\tilde\gamma_{\bm{u}}(b)]$
выпуклый и лежит выше (соответственно ниже) отрезка $[\tilde\gamma_{\bm{u}}(a)\tilde\gamma_{\bm{u}}(b)]$.
 
\parbf{Тёмная и светлая стороны.}
Пусть $\Sigma\subset\RR^3$ выпуклая поверхность,
$p\in\Sigma$ и $z\ne p$.

Точка $p$ лежит на \emph{тёмной}  (\emph{светлой}) стороне $\Sigma$ от $z$ 
если ни одна из точек $p+t\cdot(p-z)$ не лежит внутри  $\Sigma$ при $t>0$ (соответственно при $t<0$).
Пересечение тёмной стороны со светлой называется горизонтом точки $z$ и обозначается $\omega_z$.

Заметим, что если $z$ лежит внутри $\Sigma$, 
то все точки $\Sigma$ лежат на тёмной стороне от $z$ и горизонт $\omega_z$ пуст.

Если $\Sigma$ гладкая, мы можем определить внешнюю нормаль $\nu_p$ к $\Sigma$ в $p$.
В этом случае $p$ лежит на тёмной (светлой) стороне $\Sigma$ от $z$
тогда и только тогда, когда $\langle p-z,\nu_p\rangle\ge 0$
(соответственно $\langle p-z,\nu_p\rangle\le 0$).
Если  дополнительно $\Sigma$ замкнута и сильно выпукла, тогда горизонт любой точки вне $\Sigma$
образует гладкую замкнутую кривую на $\Sigma$.

Так же мы определяем светлую/тёмную сторону и горизонт в предельном случае,
когда $p$ убегает по лучу в направлении $(-\bm{u})$.

Последнее можно также определить напрямую следующим образом.
Точка $p\in\Sigma$ лежит на тёмной (светлой) стороне для вектора $\bm{u}$, 
если ни одна из точек $p+\bm{u}\cdot t$ не лежит внутри $\Sigma$ при $t>0$, (соответственно $t<0$).
Так же, как и раньше, пересечение светлой стороны с тёмной называется \emph{горизонтом} для $\bm{u}$ 
и обозначается через $\omega_{\bm{u}}$.

Значит, эти условия можно записать как $\langle \nu_p,\bm{u}\rangle\ge 0$, и $\langle \nu_p,\bm{u}\rangle\le 0$ в гладком случае.
Если $\Sigma$ ещё и замкнута и сильно выпукла, то $\omega_{\bm{u}}$ является замкнутой гладкой кривой.

\begin{thm}{Лемма Либерамна}\label{lem:liberman}
Пусть $\gamma$ есть геодезическая на выпуклой поверхности $\Sigma\subset \RR^3$.
Тогда для любой точки $z\notin\Sigma$ развёртка $\tilde\gamma_z$ локально выпукла (локально вогнута), если $\gamma$ лежит на тёмной (соответственно на светлой) стороне $\Sigma$ от $z$.

Аналогично для единичного вектора $\bm{u}$,
развёртка $\tilde\gamma_{\bm{u}}$ локально выпукла (локально вогнута), если $\gamma$ лежит на тёмной (соответственно на светлой) стороне $\Sigma$ для $\bm{u}$.
\end{thm}

Заметим, что для любой пространственной кривой $\gamma$, неравенство 
\[\tc\tilde\gamma_{\bm{u}}\le\tc\gamma.\]
выполняется для любого единичного вектора $\bm{u}$. 
С другой стороны, вариации поворотов нескольких развёрток дают оценку на вариацию поворота изначальной кривой.
Например, если $\bm{i},\bm{j},\bm{k}$ ортонормированный базис, тогда
\[\tc\gamma
\le
\tc\tilde\gamma_{\bm{i}}+\tc\tilde\gamma_{\bm{j}}+\tc\tilde\gamma_{\bm{k}}.\]

Если $\gamma$ лежит полностью на тёмной (или светлой) стороне для направления $\bm{u}$,
тогда по лемме Либермана получаем 
\[\tc\tilde\gamma_{\bm{u}}\le \pi.\]
Отсюда следует, что если $\gamma$ пересекает горизонты $\omega_{\bm{i}}$, $\omega_{\bm{j}}$ и $\omega_{\bm{k}}$
не более чем $N$ раз, тогда 
\begin{align*}
\tc\gamma
&\le
\tc\tilde\gamma_{\bm{i}}
+\tc\tilde\gamma_{\bm{j}}
+\tc\tilde\gamma_{\bm{k}}
\le
\\
&\le(N+1)\cdot \pi.
\end{align*}
Значит, если бы Основная Теорема не выполнялась для кратчайшей 
$\gamma$, то она должна была бы пересечь горизонты $\omega_{\bm{i}}$, $\omega_{\bm{j}}$ и $\omega_{\bm{k}}$ огромное число раз.

\section{Кривизна развёртки}\label{sec:curv-develop}

Пусть $\gamma\:[0,\ell]\to \Sigma$ геодезическая, параметризованная длиной
на 
замкнутой гладкой сильно выпуклой поверхности $\Sigma\subset\RR^3$.
Предположим, что для некоторого единичного вектора $\bm{u}$,
геодезическая $\gamma$ трансверсально проходит сквозь горизонт $\omega_{\bm{u}}$ в моменты времени 
$t_0<\dots <t_k$.
Определим $\alpha_i=\measuredangle(\dot\gamma(t_i),\bm{u})-\tfrac\pi2$ для каждого $i$.
Заметим, что $|\alpha_i|\le\tfrac\pi2$.

Значения $t_i$ и $\alpha_i$ 
будут называться соответственно \emph{моментами встреч} 
и \emph{углами встреч}
геодезической $\gamma$ с горизонтом $\omega_{\bm{u}}$.

Введём новое обозначение
\[\tc_{\bm{u}}\gamma\df\tc\tilde\gamma_{\bm{u}}.\]
Из леммы Либермана \ref{lem:liberman} получаем следующее.

\begin{thm}{Следствие}\label{cor:liberman}
Пусть $\Sigma\subset\RR^3$ гладкая сильно выпуклая поверхность,
$\gamma\:[0,\ell]\to \Sigma$ геодезическая, параметризованная длиной,
и $\bm{u}$ единичный вектор.
Предположим, что  $\gamma$ трансверсально пересекает горизонт $\omega_{\bm{u}}$,
а $t_0\zz<\dots <t_k$ и $\alpha_0,\dots,\alpha_k$ их моменты и углы встреч.
Тогда
\[\tc_{\bm{u}}\gamma
\le 3\cdot\pi
+
2\cdot\left|\alpha_0-\alpha_1
+\dots +(-1)^k\cdot\alpha_k\right|.
\]

\end{thm}

Далее, для получения верхней оценки на вариацию поворота,
мы получим оценку на сумму
\[\left|\alpha_0-\alpha_1
+\dots +(-1)^k\cdot\alpha_k\right|.\]
Нахождение такой оценки является основным моментом доказательства.

\parit{Доказательство.}
По лемме Либермана,
\[\tc_{\bm{u}}(\gamma|_{[t_{i-1},t_i]})=\pm(\alpha_{i-1}-\alpha_i),\]
здесь знак «$+$» 
если $\gamma_{[t_i,t_{i+1}]}$ лежит на тёмной стороне для  $\bm{u}$
и «$-$» если на светлой.
Суммируя, получаем
\[\tc_{\bm{u}}(\gamma|_{[t_0,t_{k}]})
=
\left|\alpha_0
-2\cdot\alpha_1
+\dots+(-1)^{k-1}\cdot2\cdot\alpha_{k-1}+(-1)^k\cdot\alpha_k\right|.\]
По лемме Либермана,
\[\tc_{\bm{u}}(\gamma|_{[0,t_0]}),\tc_{\bm{u}}(\gamma|_{[t_k,\ell]})\le \pi \]
Поскольку $\alpha_0,\alpha_k\le\tfrac\pi2$, утверждение следует.
\qeds 

Обозначим через $K_p$ Гауссову кривизну поверхности $\Sigma$ в её точке $p$.

\begin{wrapfigure}{r}{37 mm}
\begin{lpic}[t(-2 mm),b(-3 mm),r(0 mm),l(0 mm)]{pics/tongue(1)}
\lbl[t]{34,12;{$\omega_{\bm{u}}$}}
\lbl[lb]{26,24;{$\gamma$}}
\lbl[t]{1,11;{\tiny $\gamma(a)$}}
\lbl[t]{15,11;{\tiny $\gamma(b)$}}
\end{lpic}
\end{wrapfigure}
 
Предположим $a<b$ моменты встреч кратчайшей $\gamma$ с $\omega_{\bm{u}}$.
Дуга $\gamma|_{[a,b]}$ будет называться \emph{$\omega_{\bm{u}}$-языком},
если существует погруженный диск
$\iota\:\DD\looparrowright\Sigma$ 
такой, что замкнутая кривая $\iota|_{\partial D}$ 
образована дугами $\gamma|_{[a,b]}$ и $\omega_{\bm{u}}$.
В этом случае погружение $\iota$ называется \emph{диском языка}.

\begin{thm}{Лемма о языке}\label{lem:tongue}
Пусть $\bm{u}$ единичный вектор,
$\gamma\:[a,b]\to\Sigma$ кратчайшая, 
а также $\omega_{\bm{u}}$-язык с диском
$\iota\:\DD\looparrowright\Sigma$
на сильно выпуклой поверхности $\Sigma\subset\RR^3$.

Тогда
\[\int\limits_{\DD} K_{\iota(x)}\cdot d_{\iota(x)}\area_\Sigma\] 
принимает одно из значений
\[\alpha-\beta,
\  -\alpha+\beta,
\ \pi-\alpha-\beta,
 \pi+\alpha+\beta \pmod{2\cdot\pi},
\]
где $\alpha$ и $\beta$ обозначают углы встреч при $a$ и $b$ соответственно.

В частности
\[\bigl|\alpha-\beta\bigr|
\le
\int\limits_{\DD} K_{\iota(x)}\cdot d_{\iota(x)}\area_\Sigma.
\eqlbl{eq:key2}\]
Если дополнительно потребовать, что образ $\iota(\DD)$ 
полностью лежит на тёмной стороне для $\bm{u}$, тогда
\[\tc_{\bm{u}}\gamma
\le 
\int\limits_{\iota(\DD)} K_p\cdot d_p\area_\Sigma<2\cdot \pi.\eqlbl{eq:key3}\]
\end{thm}

\parit{Доказательство.}
Поскольку $\gamma$ геодезическая, 
параллельный перенос вдоль $\gamma$ 
переводит $\dot\gamma(a)$ в $\dot\gamma(b)$.

Заметим также, что $\bm{u}$ лежит в касательной плоскости $\Sigma$ в любой точке горизонта $\omega_{\bm{u}}$.
В частности $\bm{u}$ образует параллельное поле на $\omega_{\bm{u}}$.

Отсюда следует что параллельное перенесение вокруг $\iota|_{\partial\DD}$ 
поворачивает касательную плоскость на угол 
\[\pm(\tfrac\pi2+\alpha)\pm(\tfrac\pi2+\beta)\]
Для доказательства основного утверждения леммы остаётся применить формулу Гаусса --- Бонне.

Обозначим через $R$ правую часть \ref{eq:key2}.
Заметим, что $R\ge 0$ 
и $|\alpha|,|\beta|\le \tfrac\pi2$.
Из основного утверждения леммы следует, что  
\[R\ge\bigl|\alpha-\beta\bigr|.\]

Для доказательства \ref{eq:key3}, заметим, что в этом случае $\iota$ является вложением.
Далее заметим, что сферический образ тёмной стороны $\Sigma$ образует полусферу.
Таким образом $2\cdot\pi$ есть интеграл Гауссовой кривизны по всей тёмной стороне.
Отсюда
\[\int\limits_{\DD} K_{\iota(x)}\cdot d_{\iota(x)}\area_\Sigma
=
\int\limits_{\iota(\DD)} K_p\cdot d_p\area_\Sigma
<2\cdot\pi.\]
Применение леммы Либермана завершает доказательство.
\qeds

\section{Почти прямые кривые}

Пусть $\eps>0$.
Кривая $\gamma\:[a,b]\to\RR^3$ называется \emph{$\eps$-прямой}
если
\[(1-\eps)\cdot\length \gamma\le |\gamma(b)-\gamma(a)|.\]

\begin{thm}{Лемма}\label{lem:eps-straight}
Пусть $\eps>0$ 
и $n$ положительное целое такое, что 
$n\cdot\eps>2$.
Тогда любую кратчайшую на выпуклой поверхности $\Sigma$ в $\RR^3$
можно разбить на $\eps$-прямые дуги $\gamma_1,\dots,\gamma_n$.
\end{thm}

\parit{Доказательство.}
Пусть $\theta\in(0,\pi)$
есть число определённое равенством
\[1-\cos\theta=\eps.\]

Пусть $p,q\in \Sigma$ и $\nu_p$ и $\nu_q$ --- внешние нормали к $\Sigma$ в точках $p$ и $q$ соответственно.
Заметим, что если 
\[\measuredangle(\nu_p,\nu_q)\le 2\cdot\theta,\]
то любая кратчайшая от $p$ до $q$ в $\Sigma$
является $\eps$-прямой.

Пусть $\gamma\:[0,\ell]\to \Sigma$ кратчайшая, параметризованная длиной.

Предположим, $\gamma_{[t,\ell]}$ не $\eps$-прямая,
тогда обозначим через $t'$ максимальное число в $[t,\ell)$ такое, что дуга $\gamma_{[t,t']}$ является $\eps$-прямой.
Рассмотрим последовательность $0=t_0<t_1<\dots<t_n<\ell$ такую, что $t_{i+1}=t_i'$ для каждого $i$.
Обозначим через $\nu_i$ единичную внешнюю нормаль к $\Sigma$ в $\gamma(t_i)$. 
Из вышесказанного следует, что
\[\measuredangle(\nu_i,\nu_j)\ge2\cdot\theta\] для всех $i\ne j$.
Другими словами, открытые шары $\mathrm{B}_\theta(\nu_i)$ в сфере 
$\SS^2$ не пересекаются.

Остаётся заметить, что  
\[\area[\mathrm{B}_\theta(\nu_i)] =2\cdot\pi\cdot\eps
 \ \ \text{и}\ \ 
 \area\SS^2=4\cdot\pi.
\]
Отсюда результат.
\qeds

\begin{thm}{Следствие}\label{cor:length-diameter}
Для любой кратчайшей $\gamma$ длины $\ell$ на выпуклой поверхности $\Sigma$ в $\RR^3$ выполняется неравенство 
\[\diam \gamma\ge \tfrac\ell{10}.\]

\end{thm}

\parit{Доказательство.}
Применим лемму~\ref{lem:eps-straight} для $\eps=\tfrac12$.
\qeds

\section{Дуга в почти одном направлении}

\begin{thm}{Предложение}\label{prop:almost-const}
Для любого $\eps>0$ существует $\delta>0$ такое, что выполняется следующее условие.

Если $\gamma\:[a,b]\to\Sigma$ кратчайшая 
на гладкой сильно выпуклой поверхности $\Sigma$ в $\RR^3$,
то существует подынтервал  $[a',b']\subset[a,b]$
такой, что
\[\tc(\gamma|_{[a',b']})>\delta\cdot\tc\gamma.\]
и
\[\measuredangle(\dot\gamma(t),\bm{u})<\eps\] 
для некоторого единичного вектора $\bm{u}$
и всех $t\in[a',b']$.

Более того, при $\eps=\tfrac1{10}$ можно взять $\delta=\tfrac1{100^{100}}$.
\end{thm}

Предпошлём доказательству пару лемм.

\begin{thm}{Лемма}\label{lem:almost-const-angles}
Для любого $\eps>0$ существует $\delta>0$ такое, что выполняется следующее условие.

Предположим $\gamma$ есть кривая,
$\bm{v}_1$ и $\bm{v}_2$ два ненулевых вектора в $\RR^3$
и $0\le\alpha_1,\alpha_2\le\pi$ такие, что
\begin{align*}
\eps
&<\measuredangle(\bm{v}_1,\bm{v}_2)<\pi-\eps
\\
\alpha_i-\delta
&<\measuredangle(\bm{v}_i,\dot\gamma(t))< \alpha_i+\delta,
\end{align*}
тогда существует вектор $\bm{u}$ такой, что
$\measuredangle(\bm{u},\dot\gamma(t))<\eps$.

Более того, при $\eps<\tfrac1{10}$ 
можно взять $\delta=\eps^{10}$.
\end{thm}

Доказательство даётся прямым вычислением и мы его опускаем.

\begin{thm}{Лемма}\label{lem:almost-const}
Для любого $\eps>0$ существует $\delta>0$ такое, что выполняется следующее условие.

Предположим $\gamma\:[a,b]\to\Sigma$ является $\delta$-прямой кратчайшей на гладкой сильно выпуклой поверхности $\Sigma$ in $\RR^3$.
Пусть $\bm{v}_\gamma=\gamma(b)-\gamma(a)$.
Тогда существует подынтервал $[a',b']$ в $[a,b]$
такой, что
\[\tc(\gamma|_{[a',b']})\ge\delta\cdot\tc\gamma.\]
и
\[\alpha-\eps
\le
\measuredangle(\dot\gamma(t),\bm{v}_\gamma)
\le\alpha+\eps\] 
для некоторого фиксированного $\alpha$
и
любого $t\in[a',b']$.

Более того, при $\eps<\tfrac1{10}$ можно взять $\delta=\eps^{10}$.
\end{thm}

\parit{Доказательство.}
Без ущерба для общности, можно предположить, что
$a\zz=0$, $b=2$ и 
\[\tc(\gamma|_{[1,2]})
\ge
\tfrac12\cdot\tc\gamma.\]

Пусть $p=\gamma(0)$
и $\theta\in(0,\pi)$ такое, что $1-\cos\theta=\delta$.
Заметим, что
$$\measuredangle(\bm{v}_\gamma,\gamma(t)-p)
\le
\measuredangle(\tilde\gamma_p(1)-\tilde p,\tilde\gamma_p(2)-\tilde p)
<
\theta
\eqlbl{eq:angle1}$$
для любого $t\ge 1$.

По лемме Либермана
\[\tc_p(\gamma|_{[1,2]})< \pi+\theta.\]
Предположим $N=\lceil\tfrac\pi\theta+1\rceil$.
Тогда можно подразбить $\gamma|_{[1,2]}$ на $N$ дуг 
 $\gamma_1,\gamma_2,\zz\dots,\gamma_N$ таких, что
\[\tc_p(\gamma_n)\le \theta\eqlbl{eq:tcp}\]
для каждого $n$.

Из \ref{eq:angle1} и \ref{eq:tcp},
следует, что для каждого $n$, существует угол $\alpha_n$, такой что
\[\alpha_n-\theta
\le
\measuredangle(\dot\gamma_n(t),\bm{v}_\gamma)
\le
\alpha_n+\theta.\] 
Для доказательства основного утверждения остаётся выбрать $\delta$ так, чтобы $\theta(\delta)<\tfrac\eps{100}$
и выбрать дугу $\gamma_n$ с максимальной вариацией поворота.
\qeds

\parit{Доказательство предложения~\ref{prop:almost-const}.}
Положим $\gamma_0=\gamma$.

Для $\delta>0$, положим $n=\lceil\tfrac2\delta\rceil$.
По лемме~\ref{lem:eps-straight}, геодезическую $\gamma_0$ можно подразбить на $n$ дуг, каждая из которых $\delta$-прямая.
Выберем из них дугу $\gamma'_0$ с максимальной вариацией поворота.
Предполагая $\delta<\tfrac1{10}$, получаем
\[\tc\gamma'_0\ge\tfrac\delta{10}\cdot\tc\gamma_0.\]

Пусть $\alpha_1$ есть угол 
и $\gamma_1$ --- дуга в $\gamma'_0$, 
предоставленные леммой~\ref{lem:almost-const}.
В частности,
\begin{align*}
\tc\gamma_1&\ge\delta\cdot\tc\gamma'_0
\ge
\\
&\ge \tfrac{\delta^2}{10}\cdot\tc\gamma_0.
\end{align*}

Если
$\alpha_1\le \tfrac\eps2$ (или $\alpha_1\ge\pi-\tfrac\eps2$)
и $\delta$ достаточно мало,
тогда утверждение выполняется для дуги $\gamma_1$ и вектора $\bm{u}=\pm\bm{v}_{\gamma'_0}$.

В противном случае повторим построение для $\gamma_1$.
А именно, применим лемму~\ref{lem:eps-straight} к геодезической $\gamma_1$ и
обозначим через $\gamma_1'$ ту $\delta$-прямую дугу у которой максимальна вариация поворота.
Если значение $\delta$ мало, мы получаем 
\[
\tfrac\eps3
<
\measuredangle(\bm{v}_{\gamma_1'},\bm{v}_{\gamma'_0})
<
\pi-\tfrac\eps3
\eqlbl{eq:v1-perp-v_2}\]
Таким образом
\[\tc\gamma_1'\ge\tfrac\delta{10}\cdot\tc\gamma_1\ge \tfrac{\delta^3}{100}\cdot\tc\gamma_0.\]

Далее применим лемму~\ref{lem:almost-const} к $\gamma_1'$.
Обозначим через $\gamma_2$ и $\alpha_2$ полученные угол и дугу из $\gamma_1'$.
Снова
\[\tc\gamma_2\ge\tfrac{\delta^4}{100}\cdot\tc\gamma_0.\]

Кривая $\gamma_2$ бежит под почти постоянным углом к $\bm{v}_{\gamma'_0}$ и $\bm{v}_{\gamma'_1}$.
Неравенство \ref{eq:v1-perp-v_2}
позволяет применить лемму~\ref{lem:almost-const-angles}.
Отсюда следует основная часть.

Последняя часть доказывается прямым вычислением.
\qeds

\section{Дрейфующие геодезические}\label{sec:drifting}

В этом разделе мы вводим обозначения, используемые далее.

Выберем $(x,y,z)$-координаты в евклидовом пространстве;
обозначим через $(\bm{i},\bm{j},\bm{k})$ стандартный базис.

Плоскость параллельная, скажем $(y,z)$-координатной плоскости, будет называться $(y,z)$-плоскость.

\begin{thm}{Определение}
Гладкая кривая $\gamma$ 
называется \emph{$\bm{i}$-дрейфующей}, если оба её конца лежат на оси $x$ и 
$\measuredangle(\dot\gamma(t),\bm{i})<\tfrac1{10}$ 
для любого $t$.
\end{thm}

\parbf{$\bm{(\lambda,\mu,\nu)}$-репер.}
Пусть $\Sigma$ --- гладкая выпуклая поверхность 
и $\gamma\:[0,\ell]\to \Sigma$ --- $\bm{i}$-дрейфующая кратчайшая, 
параметризованная длиной.

Для $t\in [0,\ell]$, 
рассмотрим ориентированный ортонормированный репер $\lambda(t),\mu(t),\nu(t)$ 
такой, что вектор $\nu(t)$ есть внешняя нормаль к $\Sigma$ в $\gamma(t)$,
вектор $\mu(t)$ лежит в $(y,z)$-плоскости 
и значит вектор $\lambda(t)$ лежит в плоскости параллельной $\nu(t)$ и оси $x$.
Мы дополнительно предполагаем, что $\langle\lambda,\bm{i}\rangle\ge 0$.

Поскольку $\langle\dot\gamma(t),\bm{i}\rangle>0$, имеем 
$\nu(t)\ne\bm{i}$. 
Значит репер $(\lambda,\mu,\nu)$ определён однозначно для всех $t\in[0,\ell]$.

\parbf{Угловые функции.}
Положим
\begin{align*}
\phi(t)&=\measuredangle(\bm{i},\dot\gamma(t)),&
\psi(t)&=\tfrac\pi2-\measuredangle(\bm{i},\nu(t)),&
\theta(t)&=\tfrac\pi2-\measuredangle(\mu(t),\dot\gamma(t)).
\end{align*}
Заметим, что $|\theta(t)|,|\psi(t)|\le \tfrac\pi2$ 
и для любого $t$ существует прямоугольный сферический треугольник с катетами $|\theta(t)|,|\psi(t)|$ и гипотенузой $\phi(t)$.
В частности $\cos\theta\cdot\cos\psi=\cos\phi$.
Отсюда получаем следующее.

\begin{thm}{Утверждение}\label{clm:alpha-phi-psi}
Для любого $t$ выполняются неравенства 
\[
\phi(t)\ge |\psi(t)|\ \ \text{и}\ \ \phi(t)\ge |\theta(t)|
\]

\end{thm}

Применив лемму Либермана в направлении $\bm{i}$ получаем следующее.

\begin{thm}{Утверждение}\label{clm:7.3}
Если дуга $\gamma|_{[a,b]}$ лежит на тёмной (светлой) стороне для $\bm{i}$, 
то угловая функция $\phi$ не убывает (соответственно не возрастает) на интервале $[a,b]$.
\end{thm}

\section{Сечения}

Пусть $\gamma$ --- кривая на гладкой сильно выпуклой поверхности $\Sigma$ в $\RR^3$.
Рассмотрим плоскость $L$, 
проходящую через две точки $\gamma$, 
скажем $p=\gamma(a)$ и $q=\gamma(b)$, где $a<b$.
Пусть $L_\pm$ --- полуплоскости в $L$ ограниченные прямой, проходящей через $p$ и $q$.
Положим $\sigma_\pm=\Sigma\cap L_\pm$.

\begin{thm}{Наблюдение}\label{obs:cut}
Если $\gamma$ кратчайшая на гладкой сильно выпуклой поверхности $\Sigma\subset \RR^3$, и $a$, $b$, $\sigma_\pm$ как выше, то
\[\length\sigma_\pm\ge\length(\gamma|_{[a,b]}).\]
\end{thm}

Для доказательства достаточно заметить, что $\sigma_{+}$ и $\sigma_{-}$  суть гладкие выпуклые плоские кривые, 
соединяющие $p$ с $q$ в $\Sigma$.

Из этого наблюдения мы выведем пару оценок для дрейфующих кратчайших.

Пусть $\gamma\:[a,b]\to\RR^3$ --- кривая и $\ell$ --- прямая, не проходящая через точки на $\gamma$.
Рассмотрим непрерывную функцию $\phi\:[a,b]\to\RR$ равную азимуту $\gamma$ в цилиндрических координатах с осью $\ell$.
Если
\[|\phi(b)-\phi(a)|\ge 2\cdot n\cdot\pi,\]
мы говорим, что  $\gamma$ \emph{обошла вокруг} прямой $\ell$ хотя бы $n$ раз.

\begin{thm}{Предложение}\label{prop:around-once}
Пусть $\gamma\:[0,\ell]\to \Sigma$ есть $\bm{i}$-дрейфующая кратчайшая на выпуклой поверхности $\Sigma\subset \RR^3$. Предположим, для подынтервала $[a,b]$ в $[0,\ell]$ выполняются следующие условия:
\begin{enumerate}[(i)]
\item точки $\gamma(a)$ и $\gamma(b)$ лежат в полуплоскости ограниченной осью $x$ 
и дуга $\gamma|_{[a,b]}$ обходит вокруг оси $x$ хотя бы раз.
\item  у точки $\gamma(a)$, $x$-координата больше чем у точки $\tfrac12\cdot(\gamma(0)+\gamma(\ell))$.
\end{enumerate}
Тогда  $\gamma(b)$ лежит на тёмной стороне для $\bm{i}$.
\end{thm}

\parit{Доказательство.} 
Применим наблюдение~\ref{obs:cut} к плоскости, содержащей ось $x$, точки $\gamma(a)$ и $\gamma(b)$.

Можно предположить, что $\gamma(0)$ совпадает с началом координат 
и обе точки $p=\gamma(a)$ и $q=\gamma(b)$ лежат на $(x,z)$-координатной полуплоскости с $x\ge 0$, далее обозначаемой $\Pi$.
Мы можем также предположить, что $\sigma_+\subset \Pi$.
Обозначим через $(x_p,0,z_p)$ и $(x_q,0,z_q)$ координаты $p$ и $q$.

По условиям, $x_p<x_q<2\cdot x_p$.

Рассуждая от противного, предположим, что $\gamma(b)$ лежит на светлой стороне для $\bm{i}$.
Тогда из выпуклости кривой $\Pi\cap \Sigma$
следует, что 
\[\length\sigma_+\le \sqrt{(x_q-x_p)^2+z_p^2}.\]

С другой стороны, поскольку $\gamma|_{[a,b]}$ обходит вокруг оси $x$ хотя бы раз,
получаем 
\[\length\gamma|_{[a,b]}\ge \sqrt{(x_q-x_p)^2+(z_p+z_q)^2}.\]

Эти две оценки противоречат наблюдению~\ref{obs:cut}.
\qeds

\begin{thm}{Следствие}\label{cor:around-twice}
Если $\Sigma$, $\gamma$, $\ell$, $a$ и $b$ как в предложении~\ref{prop:around-once} 
и дуга $\gamma|_{[a,b]}$ обходит вокруг оси $x$ хотя бы дважды,
то дуга $\gamma|_{[b,\ell]}$ лежит на тёмной стороне для $\bm{i}$.
\end{thm}

\parit{Доказательство.} 
Заметим, что для любого $b'\in [b,\ell]$
существует $a'\in [a,b]$ 
такое, что условия Предложения~\ref{prop:around-once} выполняются для интервала $[a',b']$. 
Остаётся применить это предложение.
\qeds

\begin{thm}{Предложение}\label{prop:phi-psi}
Пусть $\gamma\:[0,\ell]\to \Sigma$ есть $\bm{i}$-дрейфующая кратчайшая на выпуклой поверхности $\Sigma\subset \RR^3$.
Предположим, что дуга $\gamma|_{[b,\ell]}$ 
лежит на тёмной стороне $\Sigma$ для $\bm{i}$.
Если $b\le s<t\le \ell$ и точка $\gamma(s)$ лежит в плоскости $\Pi$, проведённой через $\gamma(t)$ 
в направлении $\nu(t)$ и $\lambda(t)$,
тогда
\[\phi(s)\le \psi(t).\]
\end{thm}

\parit{Доказательство.}
Применим наблюдение~\ref{obs:cut} к плоскости $\Pi$, $p=\gamma(s)$ и $q\zz=\gamma(t)$.

Пусть $x_p$ и $x_q$ обозначают $x$-координаты $p$ и $q$.

Поскольку $\gamma|_{[s,t]}$ лежит на тёмной стороне,
её развёртка $\tilde\gamma|_{[s,t]}$ 
по отношению к $\bm{i}$ выпуклая.
В частности
\[\length(\gamma|_{[s,t]})
=
\length(\tilde\gamma|_{[s,t]})
\ge
\tfrac{x_q-x_p}{\cos\phi(s)}.\]

С другой стороны, выпуклость $\sigma_+$ влечёт, что
\[\length\sigma_+\le \tfrac{x_q-x_p}{\cos\psi(t)}.\]
Остаётся применить наблюдение~\ref{obs:cut}.
\qeds

\section{\textit{s}-пары}

Пусть  $\gamma\:[0,\ell]\to\Sigma$ --- $\bm{i}$-дрейфующая кратчайшая
на сильно выпуклой поверхности $\Sigma\subset \RR^3$.

Предположим, что $\gamma$ пересекает горизонт $\omega_{\bm{j}}$  трансверсально.

Обозначим через $t_0<t_1<\dots<t_k$ моменты встреч $\gamma$ с $\omega_{\bm{j}}$.
Положим
\begin{align*}
\phi_n&=\phi(t_n)
&
\psi_n&=\psi(t_n)
&
\theta_n&=\theta(t_n).
\end{align*}
Заметим, что $\theta_n=\pm\alpha_n$,
то есть $s_n\cdot\theta_n=(-1)^n\cdot \alpha_n$ 
для некоторой последовательности знаков $s_i=\pm1$.
В частности
\[\alpha_0-\alpha_1+\dots+(-1)^k\cdot\alpha_k
=
s_0\cdot\theta_0+s_1\cdot\theta_1+\dots+s_k\cdot\theta_k.\]
Если ориентация выбрана правильно, то  $s_n=+1$, 
когда $\nu_{\gamma(t)}$ двигается по часовой стрелке в $\SS^2$
в момент $t_n$
и если $s_n=-1$, тогда он двигается против часовой стрелки.    

Пара индексов $i<j$
образует \emph{$s$-пару}, 
если
\[
\sum_{n=i}^js_n=0\ \ 
\text{и}\ \ 
\sum_{n=i}^{j'}s_n>0,
\]
если $i<j'<j$.

Поменяв ``$+1$'' и ``$-1$'' в последовательности $s$ на скобки ``$($'' и ``$)$'' соответственно получаем, что $(i,j)$ образует $s$-пару
тогда и только тогда, когда $i$-ая скобка парна с $j$-ой скобкой.

Заметим, что каждый индекс $i$ появляется не более чем в одной $s$-паре и для любой $s$-пары $(i,j)$ имеем
\begin{itemize}
\item $s_i=1$; то есть, $i$-ая скобка должна быть открывающей.
 \item $s_j=-1$; то есть, $j$-ая скобка должна быть закрывающей.
\end{itemize}
В частности, для любой $s$-пары $(i,j)$ имеем
\begin{align*}
s_i\cdot\theta_i+s_j\cdot\theta_j&=\theta_i-\theta_j=
\\
&=(-1)^i\cdot\alpha_i+(-1)^j\cdot\alpha_j.
\end{align*}

\begin{wrapfigure}{r}{52 mm}
\begin{lpic}[t(-7 mm),b(1 mm),r(0 mm),l(0 mm)]{pics/s-pair(1)}
\lbl[br]{2,4;$+$}
\lbl[br]{3,14;$+$}
\lbl[br]{4,24;$+$}
\lbl[br]{6,34;$+$}
\lbl[bl]{16,34;$-$}
\lbl[br]{28,34;$+$}
\lbl[br]{33,44;$+$}
\lbl[bl]{43,44;$-$}
\lbl[bl]{45,34;$-$}
\lbl[bl]{46,24;$-$}
\lbl[bl]{47,14;$-$}
\lbl[bl]{47.5,4;$-$}
\lbl[lt]{4,2;$\gamma(t_i)$}
\lbl[rt]{46,2;$\gamma(t_j)$}
\end{lpic}
\end{wrapfigure}

\parbf{$\bm{s}$-пары и $\omega_{\bm{j}}$-языки.}
Пусть $(i,j)$ --- $s$-пара.
Заметим, что существует дуга в $\omega_{\bm{j}}$ от $\gamma(t_i)$ к $\gamma(t_j)$
с монотонной $x$-координатой.
Более того, существует $\omega_{\bm{j}}$-язык с этой дугой, включённой в граничную кривую.

Доказательство показано на рисунке.
Там изображено поднятие $\gamma$ в универсальное накрытие полоски в $\Sigma$ между $(y,z)$-плоскостями через $\gamma(t_i)$ и $\gamma(t_j)$.
Горизонтальные линии соответствуют поднятиям $\omega_{\bm{j}}$.

Положительное целое число $q$
назовём \emph{глубиной} $s$-пары $(i,j)$ 
(или $q\zz=\depth_{s}(i,j)$), 
если $q$ максимальное число для которого существует \emph{убывающая} последовательность $s$-пар начинающаяся $(i,j)$.
Например $s$-пара на рисунке имеет глубину $5$.

Точнее, глубина $(i,j)$ равна максимальному числу $q$
для которого существует последовательность $s$-пар,
$(i,j)=(i_1,j_1),(i_2,j_2),\zz\dots,(i_q,j_q)$ такая, что
\[i=i_1<\dots<i_q<j_q<\dots<j_1=j.\]

Заметим, что $s$-пары с одной и той же глубиной не перекрываются.
То есть, если $\depth(i,j)=\depth(i',j')$
для различных $s$-пар $(i,j)$ и $(i',j')$,
тогда $i\zz<j\zz<i'\zz<j'$ либо $i'<j'<i<j$.

Отсюда получаем следующее предложение.

\begin{thm}{Предложение}\label{prop:immersion}
Пусть $(i,j)$ --- $s$-пара.
Тогда дуги $\gamma|_{[t_i,t_j]}$ и некоторая дуга в $\omega_{\bm{j}}$ ограничивают погруженный диск в $\Sigma$, лежащий между $(y,z)$-плоскостями через $\gamma(t_i)$ и $\gamma(t_j)$.
Более того, максимальная кратность погружения не превосходит $\depth_{s}(i,j)$.
\end{thm}

\begin{thm}{Следствие}\label{cor:Sq}
Обозначим через $S_q$ подмножество индексов $\{1,\dots,k\}$
которые встречаются в $s$-парах глубины $q$.
Тогда
\[\left|\sum_{n\in S_q}(-1)^n\cdot\alpha_n\right|
=\left|\sum_{n\in S_q}s_n\cdot\theta_n\right|
\le 4\cdot\pi\cdot q.\]
\end{thm}

\parit{Доказательство.} Для индекса $n$ обозначим через $K_n$ интеграл Гауссовой кривизны куска поверхности $\Sigma$ с $x$-координатой меньше чем $x$-координата $\gamma(t_n)$.
Заметим, что
\[0\le K_1\le\dots\le K_k\le 4\cdot\pi.\]

По предложению~\ref{prop:immersion} и лемме о языке, для $s$-пары $(i,j)$ глубины $q$,
получаем
\[s_i\cdot\theta_i+s_j\cdot\theta_j=\theta_i-\theta_j\le q\cdot (K_j-K_i).\]

Утверждение следует, поскольку $s$-пары с одной и той же глубиной не перекрываются.
\qeds

\begin{thm}{Следствие}\label{cor:gamma-0}
Предположим, 
\[q=\max_{1\le i<j\le k}\left\{\biggl|\sum_{n=i}^js_n\biggr|\right\}\]
Тогда
\[\left|\sum_{n=1}^k s_n\cdot\theta_n\right|
\le 2\cdot q\cdot(q+\tfrac32)\cdot \pi.
\]
\end{thm}

\parit{Доказательство.}
Обозначим через $S$ множество всех индексов, которые появляются в $s$-парах.

Заметим, что глубина $s$-пар не превосходит $q$.
Таким образом
\[S=S_1\cup\dots\cup S_q.\]
Согласно следствию~\ref{cor:Sq},
\[\left|\sum_{n\in S}s_n\cdot\theta_n\right|\le 2\cdot q\cdot(q+1)\cdot\pi.\eqlbl{sum-S}\]

Обозначим через $R$ множество индексов, не входящих в $s$-пары,
то есть $R=\{1,\dots,k\}\backslash S$.

Для целого $r$, запишем  $i\in Q_r$, 
если
\[\sum_{n=1}^i s_n=r.\]
Заметим, что $Q_r\ne\emptyset$ для не более чем $q$ значений $r$
и каждое множество $Q_r$ содержит не более $2$ индексов не из $s$-пар;
то есть $Q_r\cap R$ содержит не более двух индексов для каждого $r$.

Поскольку $|a_n|<\tfrac\pi2$, получаем
\[\left|\sum_{n\in R}s_n\cdot\theta_n\right|
\le
q\cdot\pi.
\]
Последнее неравенство вместе с \ref{sum-S} влечёт нужное утверждение.
\qeds

\section{Геометрический рост}\label{sec:geometric-growth}

\begin{thm}{Утверждение}\label{clm:alpha-psi}
Предположим, $\psi(t)>\eps$ для $t\in[t_{i},t_{i+1}]$ и $s_i=s_{i+1}$.
Тогда 
\[|\theta_{i+1}-\theta_i|>\pi\cdot\sin\eps.\] 

\end{thm}

\parit{Доказательство.}
Заметим, что дуга $\gamma|_{[t_{i},t_{i+1}]}$ является языком, пусть $\iota\:\DD^2\zz\to \Sigma$ его диск.
Поскольку $\psi(t)>\eps$, 
сферический образ $\nu\circ\iota(\DD^2)$ 
образа $\iota(\DD^2)$ 
лежит в полукруге радиуса $\tfrac\pi2-\eps$ в $\SS^2$.
Заметим, что
\[K(\iota(\DD^2))
=
\area(\nu\circ\iota(\DD^2))
<
\pi\cdot(1-\sin\eps).\]
Остаётся применить лемму о языке~\ref{lem:tongue}.
\qeds

\begin{thm}{Утверждение}\label{clm:geometric-grouth}
Предположим $\gamma$ лежит на тёмной стороне для $\bm{i}$.
Тогда для любой пары индексов $j>i$,
таких, что
\[\biggl|\sum_{n=i}^{j}s_n\biggr|> 5\]
выполняется
\[\phi_j>\tfrac32\cdot\phi_i.\]

\end{thm}

\parit{Доказательство.}
Согласно утверждению~\ref{clm:7.3}, можно предположить, что 
\[\sum_{n=i}^{j}s_n= 6.\]

Пусть $j'$ --- последний индекс такой, что
\[\biggl|\sum_{n=i}^{j'} s_n\biggr|=5.\]

В частности $\gamma|_{[t_i,t_{j'}]}$ обходит вокруг оси $x$ хотябы дважды.
Значит для любого $b\ge t_j$ существует значение $a\in[t_i,t_{j'}]$,
такое, что $[a,b]$ удовлетворяет условиям в \ref{prop:phi-psi}.
В частности, $\psi(b)>\phi_i$ для любого $b>t_j$.
Применяя утверждение~\ref{clm:alpha-psi},
получаем, что $|\theta_j|>\tfrac\pi2\cdot \phi_i$ или
$|\theta_{j'}|\zz>\tfrac\pi2\cdot \phi_i$.
По утверждению~\ref{clm:7.3}, $\phi_n$ неубывающая и $\phi_n\ge |\theta_n|$ для каждого $n$.
В обоих случаях
\[\phi_j>\tfrac\pi2\cdot \phi_i.\]
Отсюда результат.
\qeds

\begin{thm}{Предложение}\label{prop:graph}
Если $\gamma$ является $\bm{i}$-дрейфующей кратчайшей на тёмной стороне для $\bm{i}$, тогда
\[\tc_{\bm{j}}\gamma\le 100\cdot\pi.\]

\end{thm}

\parit{Доказательство.}
Можно предположить, что
$\gamma$ пересекает $\bm{j}$-горизонт $\omega_{\bm{j}}$ трансверсально.
Пусть $t_0<\dots<t_k$ --- моменты встреч $\gamma$ с $\omega_{\bm{j}}$ и
$s_0,\dots,s_k$ --- знаки.

Напомним, что $S_q$ обозначает подмножество индексов $\{1,\dots,k\}$ из $s$-пар с глубиной $q$.
По следствию~\ref{cor:Sq},
\[\left|\sum_{n\in S_q}s_n\cdot\theta_n\right|\le 4\cdot q\cdot \pi.\]
В частности,
\[\left|\sum_{n\in S_1\cup\dots\cup S_5}
s_n\cdot\theta_n\right|
\le 
40\cdot\pi.\]

Введём обозначение $R=\{1,\dots,k\}\backslash (S_1\cup\dots\cup S_5)$ ---
это множество включает индексы из $s$-пар с глубиной хотя бы $6$, 
а также индексы, не входящие в $s$-пары.

По утверждению~\ref{clm:alpha-phi-psi},
\[\left|\sum_{n\in R}
s_n\cdot\theta_n\right|
\le
\sum_{n\in R}
|\theta_n|\le \sum_{n\in R}\phi_n.\] 
Для оценки последней суммы воспользуемся результатами из раздела~\ref{sec:geometric-growth}.
Сначала подразобьём $R$ на 5 подмножеств $R_1,\dots,R_5$, 
полагая $n\in R_m$ если $m\equiv n\pmod 5$.

Для $n\in R_m$ обозначим через $n'$ следующий индекс в $R_m$.
Значение $n'$ определено для всех $n\in R_m$, кроме последнего.
Согласно утверждению~\ref{clm:geometric-grouth}, 
$\phi_{n'}>\tfrac32\cdot \phi_n$.
То есть последовательность $(\phi_n)_{n\in R_m}$ растёт быстрее чем геометрическая с коэффициентом  $\tfrac32$.
Поскольку $\phi_n$ не убывает по $n$,
получаем
\[\sum_{n\in R_m}\phi_n< 3\cdot\phi_k.\]
Отсюда
\[\sum_{n\in R}\phi_n< 15\cdot\phi_k\le\tfrac{15}2\cdot\pi.\]

По следствию~\ref{cor:liberman},
\begin{align*}
\tc_{\bm{j}}\gamma
&\le 
2\cdot\pi
+
2\cdot [\alpha_0-\alpha_1+\dots+(-1)^k\cdot\alpha_k]< 
\\
&< 100\cdot \pi.
\end{align*}
\qeds

\section{Сборка доказательства}

Пусть $\gamma\:[0,\ell]\to \Sigma$ --- кратчайшая на выпуклой поверхности $\Sigma\subset \RR^3$.

По предложениям  \ref{prop:semicontinuity} и \ref{prop:convegence} можно предположить, что $\Sigma$  гладкая, сильно выпуклая и замкнутая,
а кратчайшая $\gamma$ имеет конечную длину.

Согласно предложению~\ref{prop:almost-const}, 
можно выбрать $(x,y,z)$-координатную систему 
и перейти к $\bm{i}$-дрейфующей дуге $\gamma'$
кратчайшей $\gamma$ такой, что
\[\tc\gamma'>\tfrac{1}{100^{100}}\cdot\tc\gamma.
\eqlbl{eq:1/100}
\]
Мы будем использовать обозначения раздела~\ref{sec:drifting} для $\gamma'$.

Подвернув $(y,z)$-координатную плоскость, можно добиться того, что
\[\tc\gamma'\le 10\cdot\tc_{\bm{j}}\gamma'\]
и что $\gamma'$ пересекает горизонт $\omega_{\bm{j}}$ трансверсально.

Согласно следствию \ref{cor:around-twice},
можно подразбить $\gamma'$ на три дуги: 
\begin{itemize}
\item \emph{Левую} $\gamma_-'$ полностью лежащую на светлой стороне для $\bm{i}$.
\item \emph{Среднюю} $\gamma_0'$ обходящую ось $x$ не больше чем $4$ раза.
\item \emph{Правую} $\gamma_+'$ полностью лежащую на тёмной стороне для $\bm{i}$,
\end{itemize}

Действительно, выберем дугу $\gamma'|_{[a,b]}$  
справа от $(y,z)$-плоскости через
$\tfrac12\cdot(\gamma'(0)+\gamma'(\ell))$,
которая делает два оборота вокруг оси $x$ 
и предположим, что значение $b$ наименьшее возможное.  
Заметим, что если $\gamma'(s)$ лежит на $(y,z)$-плоскости через $\tfrac12\cdot(\gamma'(0)+\gamma'(\ell))$
тогда $[s,b]\supset[a,b]$ 
и любая поддуга $[s,b]$ обходит вокруг оси $x$ не больше двух раз.

По следствию~\ref{cor:around-twice}, 
можно взять $\gamma'_+=\gamma'|_{[b,\ell]}$.
(В случае если интервал $[a,b]$ не определён, 
мы предполагаем, что дуги $\gamma_+'$ нет.)

Повторив построение для обратного направления оси $x$,
получим левую дугу $\gamma_-'$.
Оставшейся кусок обозначим $\gamma_0'$.
Заметим, что дуга $\gamma_0'$ разделена $(y,z)$-плоскостью через $\tfrac12\cdot(\gamma'(0)+\gamma'(\ell))$ на две дуги, 
каждая из которых обходит ось $x$ не более двух раз. 
Таким образом число обходов $\gamma_0'$ не превышает $4$.

Оценим вариацию поворота для каждой дуги $\gamma'_-$, $\gamma'_0$ и $\gamma'_+$ отдельно.

По предложению~\ref{prop:graph}, получаем 
\[\tc_{\bm{j}}\gamma'_+\le 100\cdot\pi.\eqlbl{+}\]
Аналогично  
\[\tc_{\bm{j}}\gamma'_-\le 100\cdot\pi.\eqlbl{-}\]

По следствию~\ref{cor:gamma-0},
\[\tc_{\bm{j}}\gamma'_0\le 100\cdot\pi.\eqlbl{0}\]
Оценки \ref{+}, \ref{-} и \ref{0} вместе  дают
\[\tc_{\bm{j}}\gamma'\le 300\cdot\pi.\]
Остаётся применить \ref{eq:1/100}.
\qeds

\section{Заключительные замечания}

Заметим, что основная теорема становится не верной, если в её формулировке заменить \emph{кратчайшую} на \emph{геодезическую}, даже если ограничить длину геодезической внутренним диаметром поверхности.

В дополнение к следствию~\ref{cor:length-diameter},
мы знаем пару забавных свойств кратчайших на выпуклых поверхностях, которые отличают их от геодезических. 
Хотя они нам не помогли в доказательстве, мы хотели бы их упомянуть. 
Первое из этих свойств, по сути обнаружил Анатолий Милка в \cite{milka-bending}. 
Мы дадим формулировку  только для гладких поверхностей, поскольку она проще.

\begin{thm}{Лемма о выпрямлении}
Пусть $\gamma\:[0,\ell]\to\Sigma$ --- кратчайшая, параметризованная длиной,
от точки $p$ до $q$ на выпуклой гладкой гиперповерхности $\Sigma\subset\RR^m$ и $p_t=\gamma(t)-\dot\gamma(t)\cdot t$.
Тогда, для любого $t\in [0,\ell]$, точка $q$ лежит на тёмной стороне от $p_t$.
\end{thm}

\begin{thm}{Лемма о выпуклой горбушке}
Пусть $\Sigma\subset\RR^m$ --- замкнутая выпуклая гиперповерхность и
$\Pi$ --- гиперплоскость, разрезающая $\Sigma$ на две части $\Sigma_1$ и $\Sigma_2$.
Предположим, что отражение $\Sigma_1$ в плоскости $\Pi$ заключено внутри  $\Sigma$.
Тогда $\Sigma_1$ образует выпуклое множество в $\Sigma$,
то есть любая кратчайшая в $\Sigma$ с концами в $\Sigma_1$ полностью лежит в $\Sigma_1$.
\end{thm}

Доказательства обоих утверждений можно найти в \cite{petrunin-orthodox}.

\medskip

Стыдно признаться, но мы не смогли обобщить основную теорему на старшие размерности,
то есть мы не знаем ответ на следующий вопрос.

\begin{thm}{Открытый вопрос}
Верно ли, что для любого положительного целого $m$,
существует вещественная константа $C_m$,
такая, что вариация поворота произвольной кратчайшей на выпуклой замкнутой гиперповерхности в $\RR^m$ не превосходит $C_m$?
\end{thm}

\Addresses